\documentclass[a4paper,reqno]{amsart}
\usepackage[margin=25mm,headsep=7.5mm,footskip=7.5mm]{geometry}

\usepackage{amscd,amssymb,amsmath,amsthm,amsfonts}
\usepackage{caption}


\makeatletter
\def\section{\@startsection{section}{1}%
   \z@{1.2\linespacing\@plus\linespacing}{0.8\linespacing}%
      {\normalfont\bf\boldmath\centering}}
\def\subsection{\@startsection{subsection}{1}%
   \z@{0.8\linespacing\@plus\linespacing}{0.5\linespacing}%
      {\normalfont\bf\boldmath}}
\def\@seccntformat#1{%
  \protect\textup{\protect\@secnumfont
     \ifnum\pdfstrcmp{section}{#1}=0 \bfseries\fi
     \ifnum\pdfstrcmp{subsection}{#1}=0 \bfseries\fi
     \csname the#1\endcsname\protect\@secnumpunct}}
\makeatother

\usepackage[utf8x]{inputenc}
\usepackage[english]{babel}
\usepackage{yfonts}
\usepackage[sans]{dsfont}
\usepackage{bbm}
\usepackage{mathtools}
\usepackage{accents}
\usepackage{tensor}
\usepackage[rgb]{xcolor}

\usepackage{tikz}
\usepackage{tikz-cd}
\usetikzlibrary{calc,matrix,arrows,decorations.pathmorphing}

\usepackage{calc}
\usepackage[cal=boondox,scr=boondoxo]{mathalfa}
\usepackage{enumitem}
\usepackage{marginnote}
\usepackage{booktabs}
\usepackage{scalerel}[2016/12/29]
\usepackage{url}
\usepackage{contour}
\usepackage{ulem}
\usepackage{comment}
\usepackage{placeins}

\usepackage{hyperref}
\hypersetup{colorlinks=true,pageanchor=false,linkcolor=blue,citecolor=blue,urlcolor=blue}
\usepackage[all]{hypcap}

\allowdisplaybreaks
\setlist{itemsep=2.5pt,leftmargin=30pt,labelsep=0.5em,nosep}

\newtheoremstyle{slanted}
  {2ex}{2ex}{\slshape}{\parindent}{\scshape}{.}{0.5em}{}

\theoremstyle{slanted}
\newtheorem{theorem}{Theorem}[section]

\newtheorem{proposition}[theorem]{Proposition}
\newtheorem{lemma}[theorem]{Lemma}

\newtheorem*{theorem*}{Theorem} %%%

\newtheorem{maintheorem}{Theorem}[section]

\newtheoremstyle{upright}
  {2ex}{2ex}{\upshape}{\parindent}{\scshape}{.}{0.5em}{}

\theoremstyle{upright}
\newtheorem{definition}[theorem]{Definition}

\newtheorem{remark}[theorem]{Remark}

\numberwithin{equation}{section}
\numberwithin{figure}{section}
\numberwithin{table}{section}
\makeatletter \let\c@table\c@figure \makeatother

\makeatletter
\renewenvironment{proof}[1][\proofname]{\par \pushQED{\qed}%
\normalfont \topsep6\p@\@plus6\p@\relax \trivlist
\item\relax{\slshape#1\@addpunct{.}}\hspace\labelsep\ignorespaces}
{\popQED\endtrivlist\@endpefalse}
\makeatother

\let\textit\textsl

\newcommand{\bfe}{\mathbf{e'-1}}
\newcommand{\calC}{\mathcal{C}}
\newcommand{\calH}{\text{\bf H}}
\newcommand{\cS}{\mathcal{S}}
\newcommand{\geqs}{\geqslant}

\newcommand\arxiv[2]{\href{https://arXiv.org/abs/#1}
  {\texttt{arXiv:\allowbreak #1} #2}}
\newcommand\doi[2]{\href{https://doi.org/#1}{#2}}

\newcommand{\dmnsnl}[1]{$#1$-di\-men\-sion\-al}

\newcommand{\hndlbd}[1]{$#1$-han\-dle\-bod\-y}
\newcommand{\hndlbds}[1]{$#1$-han\-dle\-bod\-ies}
\newcommand{\mnfld}[1]{$#1$-man\-i\-fold}

\newcommand{\qvlnc}[1]{$#1$-e\-quiv\-a\-lence}

\makeatletter

\def\(#1){{\upshape(\hbox to 1.2ex{\hss#1\hss})}}

\def\rel#1{\ifmmode\text{\upshape(\rel@{#1}\upshape)}%
           \else{\upshape(\rel@{#1}\upshape)}\fi}

\def\hrel#1{\ifx#1[\hrel@[\else\hrel@@{#1}\fi}
\def\hrel@[#1]#2{{\upshape(\hyperref[E:#1]{\rel@{#2}}\upshape)}}
\def\hrel@@#1{{\upshape(\hyperref[E:#1]{\rel@{#1}}\upshape)}}

\def\rel@#1{{\gdef\n@xtsp{}{\textsl{\rel@@#1@@}\n@xtsp}}}

\def\rel@@#1{\gdef\n@xt{\rel@@}%
   \ifx#1@\gdef\n@xt{\@gobble}%
   \else\ifx#1f\n@xtsp f\gdef\n@xtsp{\kern0.2ex}%
   \else\ifx#1i\n@xtsp\kern0.1ex i\gdef\n@xtsp{\kern0.1ex}%
   \else\ifx#1j\n@xtsp\kern0.2ex j\gdef\n@xtsp{\kern0.1ex}%
   \else\ifx#1l\n@xtsp\kern0.1ex l\gdef\n@xtsp{\kern0.2ex}%
   \else\ifx#1I\n@xtsp I\gdef\n@xtsp{\kern0.2ex}%
   \else\ifx#1J\n@xtsp\kern-0.1ex J\gdef\n@xtsp{\kern0.1ex}%
   \else\ifx#1'\/$\mkern0.5mu'$\gdef\n@xtsp{}%
   \else\ifx#1-{\/\upshape-}\gdef\n@xtsp{}%
   \else\n@xtsp#1\gdef\n@xtsp{}%
   \fi\fi\fi\fi\fi\fi\fi\fi\fi\n@xt}

\def\myput{\@ifnextchar[{\@put}{\@@rput[\z@,\z@][r]}}
\def\@put[#1]{\@ifnextchar[{\@@put[#1]}{\@@@@@put[#1]}}
\def\@@put[#1][{\@ifnextchar{l}{\@@lput[#1][}{\@@@put[#1][}}
\def\@@@put[#1][{\@ifnextchar{c}{\@@cput[#1][}{\@@@@put[#1][}}
\def\@@@@put[#1][{\@ifnextchar{r}{\@@rput[#1][}{\relax}}
\def\@@@@@put[{\@ifnextchar{l}{\@@lput[\z@,\z@][}{\@@@@@@put[}}
\def\@@@@@@put[{\@ifnextchar{c}{\@@cput[\z@,\z@][}{\@@@@@@@put[}}
\def\@@@@@@@put[{\@ifnextchar{r}{\@@rput[\z@,\z@][}{\@@@@@@@@put[}}
\def\@@@@@@@@put[#1]{\@@rput[#1][r]}

\let\hm@d@\leavevmode

\long\def\@@lput[#1,#2][l]#3{\setbox0\hbox{#3}\hm@d@\raise#2\hbox to\z@{\dimen0 #1%
  \advance\dimen0-\wd0\kern\dimen0\dp0\z@\ht0\z@\wd0\z@\box0\hss}\ignorespaces}
\long\def\@@cput[#1,#2][c]#3{\setbox0\hbox{#3}\hm@d@\raise#2\hbox to\z@{\dimen0 #1%
  \advance\dimen0-.5\wd0\kern\dimen0\dp0\z@\ht0\z@\wd0\z@\box0\hss}\ignorespaces}
\long\def\@@rput[#1,#2][r]#3{\setbox0\hbox{\kern#1\raise#2\hbox{#3}}%
  \dp0\z@\ht0\z@\wd0\z@\hm@d@\box0\ignorespaces}

\def\@@bold{bold}

\def\widebar#1{\text{\setbox15\hbox{$#1$}%
  \dimen15 0.45\wd15\advance\dimen15 0.15\ht15%
  \dimen16\ht15\advance\dimen16 0.00em\advance\dimen16 0.3ex%
  \dimen17 0.65\wd15\advance\dimen17 0.05\ht15\advance\dimen17 0.1ex%
  \dimen18 0.035em\advance\dimen18 0.00ex
  \@@cput[\dimen15,\dimen16][c]{\vrule depth 0pt height \dimen18 width \dimen17}}#1}

 \def\boldbar#1{\text{\setbox15\hbox{$#1$}%
   \dimen15 0.45\wd15\advance\dimen15 0.15\ht15%
   \dimen16\ht15\advance\dimen16 0.00em\advance\dimen16 0.26ex%
   \dimen17 0.65\wd15\advance\dimen17 0.05\ht15\advance\dimen17 0.1ex%
   \dimen18 0.05em\advance\dimen18 0.00ex
   \@@cput[\dimen15,\dimen16][c]{\vrule depth 0pt height \dimen18 width \dimen17}}#1}

\makeatother

\let\bar\widebar
\let\hat\widehat
\let\tilde\widetilde

\def\varemptyset{{\text{\raise.21ex\hbox{$\not$}}\mkern.15mu\mathrm{O}\mkern.15mu}}

\let\emptyset\varemptyset

\def\cal#1{\mathcal{#1}}

\newcommand{\A}{{\cal A}}
\newcommand{\C}{{\cal C}}
\newcommand{\G}{{\cal G}}

\renewcommand{\S}{{\cal S}}
\definecolor{Acolor}{cmyk}{0, 1, 0, 0.15}
\definecolor{Icolor}{cmyk}{0, 0.35, 1, 0.15}
\definecolor{Mcolor}{cmyk}{1, 0, 0, 0.25}
\definecolor{Rcolor}{cmyk}{0.5, 0, 1, 0.25}

\makeatletter

\definecolor{@@color}{cmyk}{0, 1, 1, 0.25}
\definecolor{a@color}{cmyk}{0, 1, 0, 0.15}
\definecolor{i@color}{cmyk}{0, 0.35, 1, 0.15}
\definecolor{m@color}{cmyk}{1, 0, 0, 0.25}
\definecolor{r@color}{cmyk}{0.5, 0, 1, 0.25}

\def\note{\@ifnextchar[{\@note}{\@note[]}}

\long\def\@@note[#1]#2{\def\temp{#1}%
  \ifx\temp\@empty%
    {\color{@@color}{\boldmath$\langle$}#2{\boldmath$\rangle$}}%
  \else\if\temp A%
    {\color{a@color}{\boldmath$\langle$}\textbf{A:} #2{\boldmath$\rangle$}}%
  \else\if\temp I%
    {\color{i@color}{\boldmath$\langle$}\textbf{I:} #2{\boldmath$\rangle$}}%
  \else\if\temp M%
    {\color{m@color}{\boldmath$\langle$}\textbf{M:} #2{\boldmath$\rangle$}}%
  \else\if\temp R%
    {\color{r@color}{\boldmath$\langle$}\textbf{R:} #2{\boldmath$\rangle$}}%
  \else%
    {\color{@@color}{\boldmath$\langle$}\textbf{#1:} #2{\boldmath$\rangle$}}%
  \fi\fi\fi\fi\fi}

\def\hidenotes{\long\def\@note[##1]##2{\ignorespaces}}
\def\shownotes{\let\@note\@@note}

\shownotes
\def\post{\@ifstar{\post@l}{\post@r}}

\def\post@r{\@ifnextchar[{\post@@r}{\post@@r[]}}
\def\post@l{\@ifnextchar[{\post@@l}{\post@@l[]}}

\long\def\post@@@r[#1]#2{\def\temp{#1}%
  \ifx\temp\@empty{\postitcolor{@@color}\postit{#2}}%
  \else\if\temp A{\postitcolor{a@color}\postit{\textbf{A:} #2}}%
  \else\if\temp I{\postitcolor{i@color}\postit{\textbf{I:} #2}}%
  \else\if\temp M{\postitcolor{m@color}\postit{\textbf{M:} #2}}%
  \else\if\temp R{\postitcolor{r@color}\postit{\textbf{R:} #2}}%
  \else{\postitcolor{@@color}\postit{\textbf{#1:} #2}}%
  \fi\fi\fi\fi\fi\ignorespaces}

\long\def\post@@@l[#1]#2{\def\temp{#1}%
  \ifx\temp\@empty{\postitcolor{@@color}\postit*{#2}}%
  \else\if\temp A{\postitcolor{a@color}\postit*{\textbf{A:} #2}}%
  \else\if\temp I{\postitcolor{i@color}\postit*{\textbf{I:} #2}}%
  \else\if\temp M{\postitcolor{m@color}\postit*{\textbf{M:} #2}}%
  \else\if\temp R{\postitcolor{r@color}\postit*{\textbf{R:} #2}}%
  \else{\postitcolor{@@color}\postit*{\textbf{#1:} #2}}%
  \fi\fi\fi\fi\fi\ignorespaces}

\def\hideposts{\long\def\post@@r[##1]##2{\ignorespaces}
               \long\def\post@@l[##1]##2{\ignorespaces}}
\def\showposts{\let\post@@r\post@@@r\let\post@@l\post@@@l}

\hideposts

\def\postitfontsize#1{\let\postitf@ntsize#1}
\def\postitcolor#1{\def\postitc@lor{#1}}

\postitfontsize{\tiny}
\postitcolor{red}

\newdimen\postitm@xwidth
\def\postitmaxwidth#1{\postitm@xwidth=#1}

\postitmaxwidth{8cm}

\def\postit{\@ifstar{\postit@l}{\postit@r}}

\makeatother

\newcommand{\N}{\mathbb{N}}
\newcommand{\Z}{\mathbb{Z}}
\newcommand{\R}{\mathbb{R}}

\newcommand{\Fopp}{\text{\bf F}^\text{opp}}
\newcommand{\Free}{\text{\bf F}}

\newcommand{\diam}{\otimes}

\newcommand{\bcsum}{\mathbin{\natural}}

\newcommand{\Obj}{\mathop{\mathrm{Obj}}\nolimits}

\newcommand{\one}{{\boldsymbol{\mathbbm 1}}}
\newcommand{\idone}{\id_\one}

\newcommand{\RHB}{4\mathrm{HB}}

\newcommand{\pres}{\cal P}

\newcommand{\Algf}{4\mathrm{Alg}}

\newcommand{\calg}{2\mathrm{Alg}}

\newcommand{\falg}{\mathfrak{q}}
\newcommand{\ffree}{\mathfrak{f}}
\newcommand{\ftop}{\mathfrak{r}}
\newcommand{\fh}{\mathfrak{h}}

\newcommand{\KTf}{\mathrm{4KT}}

\newcommand{\cwcob}{2\mathrm{CW}}
\newcommand{\emptys}{{\varemptyset}}

\newcommand{\perm}[1]{{\mathfrak S}_{#1}}

\newcommand{\Deltap}{\dot{\Delta}}
\newcommand{\epsp}{\dot{\epsilon}}
\newcommand{\multp}{\dot{\m}}
\newcommand{\etap}{\dot{\eta}}
\newcommand{\lp}{\dot{\lambda}}
\newcommand{\Lp}{\dot{\Lambda}}
\newcommand{\Sp}{\dot{S}}
\newcommand{\gammap}{\dot{\braid}}
\newcommand{\idp}{\dot{\id}}

\newcommand{\beqn}{\begin{eqnarray*}}
\newcommand{\eeqn}{\end{eqnarray*}}
\newcommand{\be}{\begin{equation}}
\newcommand{\ee}{\end{equation}}
\newcommand{\ba}{\begin{array}}
\newcommand{\ea}{\end{array}}

\newcommand{\pibar}{\overline{\Pi}}
\newcommand{\m}{\mu}
\newcommand{\defeq}{\stackrel{\text{\rm def}}{=}}

\newcommand{\expanse}{\!\nearrow\!}
\newcommand{\collapse}{\!\searrow\!}
\newcommand{\deform}{\diagup\mkern-7mu\searrow}

\newenvironment{sparagraph}[1]{\noindent {\sl #1.} }{}

\newcommand{\id}{\mathrm{id}}

\newcommand{\ev}{\mathrm{ev}}
\newcommand{\coev}{\mathrm{coev}}

\newcommand{\braid}{\gamma}
\newcommand{\prodH}{\mu}
\newcommand{\unitH}{\eta}
\newcommand{\coprH}{\Delta}
\newcommand{\counH}{\varepsilon}
\newcommand{\antipH}{S}

\newcommand{\ribmorH}{\tau}
\newcommand{\copairH}{w}

\newcommand{\cadjH}{\mathcal{A}\kern-0.3em\mathcal{d}}
\newcommand{\idfunct}{\mathcal{I}\kern-0.45em\mathcal{d}}

\newcommand{\monH}{\Omega}

\usepackage{accents}
\newlength{\dhatheight}

\begin{document}

\title[
Algebraic presentation of the bordism category of 2-dimensional CW-complexes]{Algebraic presentation of the bordism category of 2-dimensional CW-complexes} 
%\label{Version 0.9(50) / \today}}
%
\author[I. Bobtcheva]{Ivelina Bobtcheva} 
\address{Institute of Mathematics, University of Zurich, Winterthurerstrasse 190, CH-8057 Zurich, Switzerland}
\email{ivelina.bobtcheva@math.uzh.ch}%
%\date{}

\maketitle

\begin{abstract}
We prove that $S^1$ carries the structure of a unimodular cocommutative Hopf algebra in the bordism category $\cwcob $ of 2-equivalence classes of  2-dimensional CW-complexes and  that $\cwcob $ is actually equivalent to the symmetric monoidal category freely generated by such a Hopf algebra.  The proof proceeds by constructing an equivalence between $\cwcob$ and a category of finite relative group presentations.
Finally, we show that the algebraic structure of the category $\RHB $ of 2-equivalence classes of orientable 4-dimensional 2-handlebodies  refines that of $\cwcob $.
\end{abstract}



\section*{Introduction%
\label{introduction/sec}}

One of the main open problems in  two-dimensional homotopy theory is  the long-standing {\sl generalized Andrews-Curtis (AC) conjecture}, which asserts that any simple homotopy equivalence between two 2-complexes can be obtained through a 2-deformation consisting of changing the attaching maps of the 2-cells by homotopy and expansions and collapses of disks of dimension at most two \footnote{Here we follow the terminology in \cite{Wa66, Qu95, BBDP23}, which is motivated by the observation  that  two $n$-dimensional CW-complexes are $n$-equivalent if one can be deformed to the other through $n$-dimensional complexes. Such terminology is not universally accepted. For example in \cite{KKN25, HMS93}, this is  called  $n+1$-deformation, since it can be achieved through expansions and collapses of disks of dimension at most $n+1$.} \cite{AC65, HMS93}. Observe that while in dimensions $n \geq 3$  two $n$-complexes are simple homotopy equivalent if and only if they are related by an $n$-deformation \cite{Wa66, Wa80}, the AC-conjecture is expected to be false and different proposals for counterexamples have been made. An extensive reference on this problem is \cite{HMS93}. What makes the AC-phenomenon difficult to detect is its instability with respect to taking one-point union with copies of $S^2$. Indeed, if two 2-complexes $X_1$ and $X_2$ are simple homotopy equivalent, then for some $k\geq 0$ $X_1\vee ^k S^2$ and $X_2\vee ^k S^2$ are 2-equivalent \cite{Wa66}. This phenomenon is similar to the one observed in the case of 4-manifolds: if two smooth oriented compact 4-manifolds $M_1$ and $M_2$ are homeomorphic, then for some $k\geq 0$, $M_1\# ^k S^2\times S^2$ and $M_2\# ^k S^2\times S^2$ are diffeomorphic \cite{Go84}. In \cite{KKN25} this analogy is further extended to the universal pairing by showing that, similarly to the case of simply-connected 4-manifolds, the universal pairing for 2-complexes is not positive, which implies that  unitary TQFTs of 2-complexes do not detect the difference between simple homotopy equivalence and 2-deformations. 

In the present work we classify all TQFTs on 2-equivalence classes of 2-complexes by giving an algebraic presentation of the corresponding bordism category $\cwcob$. TQFTs of this kind were first constructed by Frank Quinn in \cite{Qu95} and the reader can find in \cite{KK17} a summary of the extensive study of Quinn's theory. The category $\cwcob$ with which we work is a skeleton of Quinn's bordism category and  has as objects   1-complexes with a single 0-cell, while the set of  morphisms $Y_1\to Y_2$ consists of 2-equivalence classes of pairs of connected 2-complexes $(X, Y_1\vee Y_2)$. In particular, $\cwcob$ is a PROP, i.e. it has a symmetric strict monoidal structure given by one-point union of complexes and any object is the monoidal power of a single generating object: $S^1$. It is well known  that 2-equivalence classes of 2-complexes can be described by group presentations modulo a set of combinatorial moves \cite[Chapter I, Theorem 2.4]{HMS93}. We use the relative version of this result for CW-pairs in order to show that there is a category equivalence between $\cwcob$ and a category of  relative presentations $\pres$. Then the simple combinatorial nature of the morphisms in $\pres$ allows us to give the following algebraic presentation of $\cwcob$. 

\begin{maintheorem}\label{main1/thm}
Let $\calg$ be the symmetric monoidal category freely generated by a  unimodular cocommutative  Hopf algebra. Then the functor $\Psi: \calg \to \cwcob$ sending the generating Hopf algebra of $\calg$ to $S^1$  is an equivalence of symmetric monoidal categories. 
\end{maintheorem} 

In \cite{Pi02}, Pirashvili proves that the opposite of the category of finitely generated free groups $\Fopp$ is equivalent to the symmetric monoidal category $\calH$ freely generated by a cocommutative Hopf algebra. A simpler combinatorial proof of Pirashvili's result is given by Habiro \cite{H16}. Theorem \ref{main1/thm} can be viewed as an extension of Pirashvili's result to finite presentations of groups, and our approach is similar to Habiro's.

In view of the analogy between 2-complexes and smooth 4-dimensional manifolds, mentioned above, it is interesting to understand  the relationship between the algebraic structure of the bordism category of 2-complexes and that of their 4-dimensional thickenings, i.e. of the category $\RHB$ of 2-equivalence classes of 4-dimensional 2-handlebodies studied in \cite{BP11, BD21, BBDP23}. The 2-equivalence relation on handlebodies consists in diffeomorphisms implemented by a sequence of handle slides and creation/cancellation of canceling pairs of 1- and 2-handles. Similarly to the case of 2-complexes, it is expected that 2-deformations form a proper subclass of diffeomorphisms \cite{Go91}.

The category $\RHB$ is a PROB, i.e. a strict braided monoidal category in which every object is a power of a single generating object: the solid torus. Building on works of Crane, Yetter and Kerler \cite{CY94, Ke02}, a complete algebraic presentation of $\RHB$ is given in \cite{BP11} by showing that the torus admits the structure of a BP Hopf algebra, which is a braided unimodular Hopf algebra satisfying some additional ribbon axioms  (see Section \ref{Hr-Htr/sec}). Moreover, if $\Algf$ is the strict braided monoidal  category freely generated by a  BP Hopf algebra,  there is a braided  monoidal equivalence $\Phi: \Algf \to \RHB$ (see \cite{BBDP23} for a simpler direct proof of this result).

The deformation retraction of each 4-dimensional 2-handlebody   onto its spine gives a  well-defined braided monoidal functor  $\ftop:\RHB \to\cwcob$ which, through the identification of  $\cwcob$ and $\pres$, can also be seen as a forgetful functor reading  the relative presentation of the fundamental group of the handlebody from the attaching maps of its handles. Under this functor, the algebraic structure of $\RHB$ can be seen as  a refinement of the algebraic structure of $\cwcob$.

\begin{maintheorem}\label{main2/thm}
There is a braided monoidal quotient functor $\falg:\Algf\to \calg$ such that the following diagram commutes:
\[
\begin{tikzcd}
\Algf \arrow{r}{\Phi} \arrow[swap]{d}{\falg} & \RHB \arrow{d}{\ftop} \\%
\calg \arrow{r}{\Psi}& \cwcob
\end{tikzcd}
\]
\end{maintheorem} 

At the level of numerical invariants, a relation between quantum invariants  of 4-dimensional thickenings and invariants of their 2-dimensional spines was studied in \cite{BQ05}, where to any finite semisimple ribbon category is associated an invariant of 4-thickenings which,  if the category is symmetric, depends only on the spine and  second Whitney number of the thickening. Theorem \ref{main2/thm}   can be seen as the universal algebraic counterpart of this result. %and \cite[Corollary B]{BBDP23}

\smallskip
{\sl Organization.}
In Sections \ref{cw/sec} and \ref{pres-category/sec} we introduce the  category $\cwcob$ and show that it is equivalent to the category $\pres$ of  relative presentations modulo AC-moves.  Then in Section \ref{Hopfc/sec} we define a unimodular, cocommutative  Hopf algebra $H$ in a symmetric monoidal category $\C$ and introduce the category $\calg$, freely generated by such an algebra.
Sections \ref{hopf-pres/sec} and \ref{Hr-Htr/sec} are devoted respectively to the proofs of Theorems \ref{main1/thm} and \ref{main2/thm}. Finally, in Appendix \ref{faith/sec} we explain why Theorem \ref{main1/thm} can be considered as a generalization of  Pirashvili's theorem to finitely presented groups. 

\smallskip
 The paper is intended as an introduction to the topic of  algebraic presentations of bordism categories. Therefore the exposition is self-contained and accompanied by numerous illustrations. This choice is motivated by the fact that the category of 2-complexes exhibits a sufficiently simple yet rich algebraic structure in which the product and coproduct are not dual to each other. Furthermore, thanks to Theorem \ref{main2/thm}, it also represents a good starting point for the study of the algebraic presentations of (co)boundary categories in  dimensions three and four.
 
\smallskip
{\sl Acknowledgments.} The idea of a possible relationship between the quantum invariants of 3-dimensional cobordisms, 4-dimensional 2-handlebodies and 2-dimensional CW-complexes is due to Frank Quinn who started this line of research and whose insights have remained the main driving force behind it. The author would also like to thank Riccardo Piergallini for carefully reading and indicating some gaps in a preliminary version of this work, Gw\'{e}na\"{e}l Massuyeau for explaining to her the connection with Pirashvili's result and Anna Beliakova for stimulating discussions. The author is supported by the NCCR SwissMAP and Grant $200020\_207374$ of the Swiss National Science Foundation and the Simons Collaboration ''New Structures in Low-Dimensional Topology''.

%%%%%%%%%%%%%%%%%%%
%%%%%%%%%%%%%%%%%%%%

\section{Preliminaries on 2-dimensional CW-complexes%
\label{cw/sec}}

We recall some basic definitions and facts about CW-complexes. The reader can find all details and a complete reference list in \cite{HMS93}. 

\medskip
\subsection{CW complexes} An {\sl $n$-dimensional CW complex} or simply an $n$-complex, is a space $X=X^n$ constructed inductively in the following way:

(1) $X^0$ is a discrete space, each point being a 0-cell.

(2) $X^k$ is obtained by attaching to $X^{k-1}$ a disjoint family $D_i^k$ of closed $k$-disks via continuous functions $\varphi_i:\partial D_i^{k}\to X^{k-1}$. In particular, $X^k$ is the quotient space of $X^{k-1}\sqcup _i D_i^k$ under the identification $x\sim \varphi_i(x)$ for $x\in \partial D_i^k$. The homeomorphic image of the interior of $D_i^k$ under the restriction of the quotient map $\phi_i:D_i^k\to X^k$ is the {\sl $k$-cell} $e_i^k$. Then $\varphi_i$ and $\phi_i$ are called respectively the {\sl attaching map} and the {\sl characteristic map} of $e_i^k$. 

The space $X$ is endowed with the weak topology and it is considered up to {\sl cellular isomorphism}, i.e. up to homeomorphism which maps open cells to open cells. In particular, the attaching and the characteristic maps are not unique, but two attaching (resp. characteristic) maps of a given cell are homotopic up to composition with an orientation reversing homeomorphism of $\partial D$ (resp $D$) to itself.

 We will assume all CW complexes to be finite, i.e. they  consist of finitely many cells. 
\smallskip

\subsection{$n$-deformation and simple homotopy}
A {\sl subcomplex}  of a CW complex $X$ is a closed subspace $Y\subset X$ that is a union of cells of $X$. For any $0\leq k\leq n$, the subspace $X^k$ is a subcomplex of $X$, called the {\sl $k$-skeleton} of $X$. A pair $(X,Y)$ consisting of an $n$-complex $X$ and a subcomplex $Y$ will be called a pair of $n$-complexes. The pair will be called {\sl connected} if both $X$ and $Y$ are connected.

\begin{definition}\label{collapse/defn} A subcomplex $Y$ is obtained from $X$ by an {\sl elementary collapse} of dimension $k$, or equivalently $X$ is obtained from $Y$ by an {\sl elementary expansion}, if there exists a pair $(e^k,e^{k-1})$ of cells of $X$ with characteristic maps $\phi_k: D^k\to X$ and $\phi_{k-1}: D^{k-1}\to X$ and an embedding $\iota: D^{k-1}\to \partial D^k$, such that
\begin{itemize}
\item[(i)] $X=Y\sqcup e^k\sqcup e^{k-1}$ and
\item[(ii)] $ \phi_{k-1}=\phi_k\circ\iota$ while $\phi_k(\partial D^k-\iota( D^{k-1}))\subset Y$.
\end{itemize}
Here $e^{k-1}$ is called a {\sl free face} of $e^k$ and of $X$. A sequence of elementary expansions (resp. collapses) is called an {\sl expansion} (resp. a {\sl collapse}) and denoted  by $Y\expanse X$ (resp. $X\collapse Y)$, while a sequence where both expansion and collapses are allowed is called a {\sl deformation} and denoted by $X\deform Y$.
\end{definition}

\begin{definition}\label{simplehomeq/defn} We say that two  CW-pairs $(X,Y)$ and $(X',Y)$ are {\sl simple-homotopy equivalent} if there is a map $X\to X'$ which is homotopic to a deformation and induces the identity on the common subcomplex $Y$.
Two CW-complexes $X$ and $X'$ are  simple-homotopy equivalent if the CW-pairs 
$(X,\emptys)$ and $(X',\emptys)$ are  simple-homotopy equivalent.
\end{definition}

\begin{definition}\label{ndef/defn} We say that two $n$-dimensional CW-pairs $(X,Y)$ and $(X',Y)$ are {\sl $n$-equivalent} if one can be deformed to the other by homotopies of attaching maps of the cells and expansions and collapses of dimension at most $n$, which restrict to the identity on the common subcomplex $Y$. A sequence of such moves is called an {\sl $n$-deformation}. 

Two $n$-dimensional CW-complexes $X$ and $X'$ are $n$-equivalent  if the CW-pairs 
$(X,\emptys)$ and $(X',\emptys)$ are $n$-equivalent.
\end{definition}

In the rest of the paper, unless otherwise specified, we will work only with 2-dimensional CW-pairs $(X,Y)$ with $Y\subset X^1$. 

%%%%%%%%%%
\subsection{Elementary 2-deformation moves}
\label{el-moves/sec}

\begin{figure}[htb]
  \includegraphics{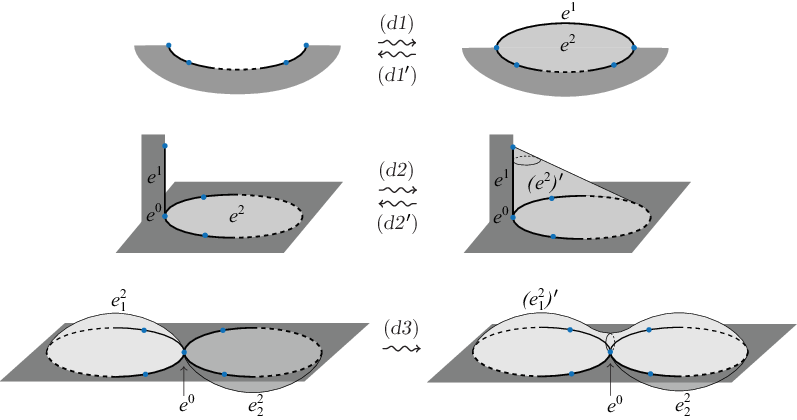}
\caption{Elementary 2-deformation moves}
\label{twodeformation/fig}
\vskip-6pt
\end{figure}

We will focus now on three particular 2-deformation moves between the CW-pairs $(X,Y)$ and $(X',Y)$ and, as it will be shown in Proposition \ref{standard/theo}, any 2-deformation between two standard 2-complexes is a composition of those moves. The moves are illustrated in Figure \ref{twodeformation/fig}, where only the 0-, 1- and 2-cells involved are shown, while the rest of the CW-complex is presented in dark-gray. The reader should keep in mind that  some of the drawn 0- and 1-cells may coincide with each other.  

\medskip
\noindent{\hrel{d1-1'}} \label{E:d1-1'}  \label{E:d1}\label{E:d1'}Expanding/collapsing a pair $(e^2,e^1)$, $e^1\notin Y$ (see Definition \ref{collapse/defn}).

\medskip
\noindent { \hrel{d2-2'}}\label{E:d2-2'} \label{E:d2}\label{E:d2'} Pushing a 2-cell $e^2$ along/off the 1-cell $e^1$.

\begin{figure}[htb]
\centerline{\includegraphics{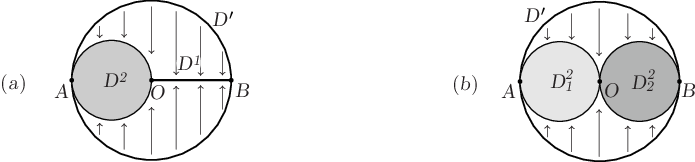}}
\caption{}
\label{slidings/fig}
\vskip-6pt
\end{figure}

\smallskip
The move is defined if the attaching maps of the 1-cell $e^1$ and of the 2-cell $e^2$ send  $P_1\in \partial D^1$ and $P_2\in \partial D^2$ to the same 0-cell $e^0$. Consider the disk $D'$ with radius 1 and center $O$ presented in Figure \ref{slidings/fig} (a). Let $E\subset D'$ be the union of the darker disk of radius $1/2$  and $OB$. We identify $E$ with $D^2\vee_{P_1=P_2} D^1$ and consider a map $\tilde \phi: E\to X$, whose restriction to $D^2$ (resp. $D^1$) coincides with the characteristic map of $e^2$ (resp. $e^1$). Let also $r: D'\to E$ be the retraction which maps each point of $D'$ to  $E$ as shown in Figure \ref{slidings/fig} (a), and $f_t:D'\to D'$ be the homotopy which maps any point $x\in D'$ to the point $f_t(x)$ dividing the segment from $A$ to $x$ in ratio $1:(1-t)$. Then $\tilde \phi\circ r\circ  f_t|_{S^1}$ changes the attaching map of $e^2$  by homotopy and the 2-complex for $t=1$ is the result of pushing $e^2$ along $e^1$.  
Observe that such a complex depends not only on the 1- and 2-cells in question, but also on the choice of particular preimages of $e^0$ under the attaching maps of such cells.

The inverse of such a move, to which we refer as {\sl pushing $e^2$ off $e^1$}, is defined when the attaching map of $e^2$ maps two consecutive arcs in $\partial D^2$ onto  $e^1\cup \bar{e}^1$, where $\bar{e}^1$ is the 1-cell $e^1$ taken with the opposite orientation.   

\medskip
\noindent { \hrel{d3}}\label{E:d3}  Sliding a 2-cell $e^2_1$ over another 2-cell $e^2_2$.

\smallskip
Suppose that a CW-complex $X$ has 2-cells $e^2_1$ and  $e^2_2$, whose closures contain a common 0-cell $e^0$ and let $P_1$ (resp. $P_2$) be  a preimage of $e^0$ under the attaching map of $e^2_1$ (resp. $e^2_2$). Then the definition of the sliding of  $e^2_1$ over  $e^2_2$ follows the steps in  \hrel{d2-2'}, where  the segment $OB$ is replaced by the disk $D_2^2$ in Figure \ref{slidings/fig} (b),  $E=D_1^2\cup D_2^2\subset D'$  is identified with $D^2\vee_{P_1=P_2} D^2$ and the map $\tilde \phi: E\to X$ is defined by requiring that its  restriction  to the first (resp. second) copy of $D^2$ coincides with the characteristic map of $e^2_1$ (resp. $e^2_2$). Observe that once again the final result of the move depends on the choice of preimages of $e^0$. Moreover, one can obtain the inverse move by sliding the resulting 2-cell once more over $e^2_2$ and then pushing it off the 1-cells in the boundary of $e^2_2$. We will discuss this later in more detail in the case of standard complexes.

\smallskip
In what follows, it would be useful to have a name for a fourth move which is obtained as  composition of the previous ones.

\smallskip
\noindent { \hrel{d4}}\label{E:d4} Sliding a 1-cell $e^1_1\notin Y$ over another 1-cell $e^1_2$.

\begin{figure}[htb]
\centerline{\includegraphics{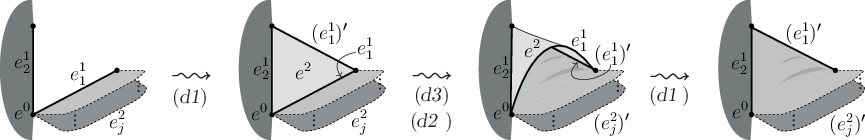}}
\caption{Definition of the 2-deformation move \hrel{d4}} % 
\label{one-sliding/fig}
\end{figure}

\smallskip
The move is defined when the images of the attaching maps of two  1-cells $e^1_i$, $i=1,2$ contain a common 0-cell $e^0$ as follows (cf. Figure \ref{one-sliding/fig}).
We first expand a pair of 2- and 1- cells $(e^2,e^1)$, where $\partial \text{Cl}( e^2)=\text{Cl}( (e_1^1)'\cup e^1_1\cup  e^1_2)$.  Then we use move \hrel{d3} to slide  all the 2-cells $e_j^2$ over $e^2$, attached to $e_1^1$, and using move \hrel{d2'} we push  them off $e_1^1$ making it free. Finally, we  collapse the pair $(e^2,e^1_1)$.

%%%%%%%%%%HERE (summer 2020)
\subsection{Standard complexes}

Standard complexes were introduced by Wright in \cite{Wr75} (see also \cite[Chapter I]{HMS93}) in order to prove the correspondence between 2-complexes and group presentations. We will need them in the slightly more general context of CW-pairs $(X,Y)$, $Y\subset X^1$.

\begin{definition}\label{standardCW/defn} A 2-dimensional CW-pair $(X,Y)$,  $Y\subset X^1$ will be called {\sl standard} if it satisfies the following properties:
\begin{itemize}
\item[(i)] $Y^0=X^0$ consists of a single 0-cell (vertex) and $X^1$ is the union of $X^0$ and a (finite) collection $e_i^1$ of 1-cells whose boundaries are attached to it;
\item[(ii)] $X^2$ is obtained by attaching to $X^1$ a (finite) collection of 2-cells $e_j^2$, where each $\partial D_j^2$ is subdivided into edges and vertices; each vertex of $\partial D_j^2$ is sent by the corresponding attaching map to $X^0$ and each open edge of $\partial D_j^2$ is sent either to $X^0$ or homeomorphically to some open $e_i^1$.
\end{itemize}
If $Y=X^0$, we will refer to the corresponding standard CW-pair $(X,X^0)$ simply as the {\sl standard 2-dimensional CW-complex} $X$. 

Given two standard CW-pairs $(X,Y)$ and $(X',Y')$, we define their one-point union $(X\vee X', Y\vee Y')$ to be the standard  CW-pair obtained by  identifying  $X^0$ and ${X'}^0$ in the disjoint union $(X,Y)\sqcup (X',Y')$. Observe that this operation defines a monoidal structure on the set of  standard  CW-pairs whose identity is given by the trivial complex consisting of a single 0-cell.
\end{definition}

%\vspace{-5mm}
\begin{remark} \label{pres/remark} To any standard 2-dimensional CW complex $X$ we can associate (in a non-unique way)  a presentation of its fundamental group in terms of generators and relations. In order to do this, we need to make a choice of an orientation of each 1-cell and a choice of a base point and an orientation of $\partial D^2$ for each 2-cell in such a way that the attaching map $\phi:\partial D^2\to X^1$ of the 2-cell maps its base point to $e^0$. Then the desired presentation has one generator for each 1-cell and one relator for each 2-cell, obtained by reading off the word corresponding to $\phi(\partial D)$ and different choices of orientations or of base points change the presentation by substitutions of  generators or relators  by their inverses and by cyclically  permuting the letters of the relations. Analogously, to any standard CW-pair $(X,Y)$ we can associate a  presentation of the fundamental group of $X$, in which the generators corresponding to $Y$ have been  marked. 
\end{remark}

\smallskip
Observe that the set of standard pairs is invariant under the moves  \hrel{d1-1'}, \hrel{d2-2'} and \hrel{d3}. Moreover, the following proposition implies that in the study of 2-equivalence classes of  CW-pairs we can limit ourselves to the study of standard pairs modulo such moves. We only sketch the proof which is a straightforward adaptation of Lemma 4 in \cite{Wr75} and Theorem 2.4 in Chapter I of \cite{HMS93} to the case of CW-pairs.

%\newpage
\begin{proposition}\label{standard/theo} Any connected  CW-pair $(X,Y)$ in which $Y\subset X^1$ is a standard one-dimensional complex, is 2-equivalent relative to $Y$ to a standard pair. Moreover, if two standard pairs are 2-equivalent, then there is a 2-deformation between them which is a composition of the moves \hrel{d1-1'}, \hrel{d2-2'} and \hrel{d3}.
\end{proposition}

% For more details see draft.rtf 2/1/25

\begin{proof} Given  a spanning tree $T$  of the 1-skeleton of $X$ with root  $v=Y^0$, for each 1-cell $e_i^1\in X^1\smallsetminus T$, $i=1,\dots n$, we  choose an orientation and two (possibly coinciding)  paths  $l_i'$, $l_i''$ in $T$ from the endpoints of $e_i^1$ to $v$, where if $e_i\in Y$, $l_i'$ and $l_i''$ are the constant paths. Then $\left\{g_i=(l_i')^{-1}\circ e_i^1\circ l_i''\right\}_{i=1}^n$ represents a set of generators of the fundamental group of $X^1$ (and of $X$) and we construct  a standard CW-pair $(X_T,Y)$, together with a 2-deformation from $(X,Y)$ to  $(X_T,Y)$, by first sliding each $e_i^1$ over all 1-cells in $l_i'$ and $l_i''$ (move \hrel{d4}), thereby obtaining a CW-complex whose 1-skeleton consists of the union of $T$,  and $n$ loops based at $v$, then using   \hrel{d2'} to push  all 2-cells off $T$ and eventually contracting $T$ to $v$ by collapses of 0/1-cells. Up to homotopy of the attaching maps of the 2-cells in the 1-skeleton in order to achieve the condition \(ii) in Definition \ref{standardCW/defn}, the resulting pair is $(X_T,Y)$.
Observe that different choices of paths and of spanning tree change $(X_T,Y)$ by a sequence of slides of 1-cells (move \hrel{d4}). 

Now if $(X',Y)$ is obtained from  $(X,Y)$ by an expansion or by a collapse of a canceling 1-2 pair (resp. 0-1 pair) by an appropriate choice of the spanning trees of $X$ and $X'$ we can obtain that  $(X_T,Y)$ and $(X'_{T'},Y)$ are related  by the expansion or the collapse of the same 1-2 pair (resp. by the identity). 
On the other hand, changing the attaching map of a given 2-cell by homotopy induces a homotopy $\varphi_t:\partial D \to X_T\smallsetminus e^2$, $0\leq t\leq 1$, of the attaching map of the corresponding 2-cell $e^2$ in $X_T$. Then $\varphi_{0}(\partial D)$  and $\varphi_{1}(\partial D)$ are related by conjugations and multiplications by elements in the normal subgroup generated by the attaching maps of the 2-cells different from $e^2$, and those transformations can be realized as a sequence of moves  \hrel{d2} and  \hrel{d3}. \end{proof}

%%%%%%%%%%
\subsection{The category $\cwcob$ of  2-dimensional CW-complexes%
\label{cw-category/sec}}

We can now define the category  of 2-dimensional CW-complexes. 

\medskip
\noindent{\sl Objects.}$\;$  The objects   of $\cwcob$ are {\sl standard} 1-dimensional complexes equipped with an  orientation and enumeration %\footnote{A finite enumerated set is a set $A$ together with a bijective function $\{1,2, \ldots, n\}\to A$.}  
of their 1-cells, considered up to homeomorphism which preserves the extra structure. We will use the notation $\bar Y_n$ to indicate the object of $\cwcob$ whose underlying  CW-complex $Y_n$ has $n$ 1-cells. 

\medskip
\noindent{\sl Morphisms.}$\;$ A morphism $(X,\bar Y_{n_1}, \bar Y_{n_2}): \bar Y_{n_1}\to  \bar Y_{n_2}$ in $\cwcob$ is a 2-equivalence class of standard 2-dimensional  CW-pairs   $(X,Y_{n_1}\vee Y_{n_2})$ in which $Y_{n_1}$ and $Y_{n_2}$ are equipped with  orientation and  enumeration of their 1-cells.  The composition of two morphisms $(X,\bar Y_{n_1}, \bar Y_{n_2})$ and $(X',\bar Y_{n_2}, \bar Y_{n_3})$ is given by $(X\sqcup_{Y_{n_2}}X',\bar Y_{n_1},\bar Y_{n_3})$, where $X\sqcup_{Y_{n_2}}X'$ is the gluing of $X$ and $X'$ along $Y_{n_2}$ by a homeomorphism which respects the enumeration and the orientation of the 1-cells.

\begin{figure}[hbt]
  \includegraphics{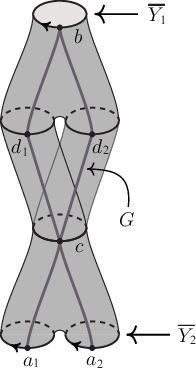}
\caption{A morphism $(X, \overline{Y}_2, \overline{Y}_1)$ in  $\cwcob$}
\label{example-cwcob/fig}
\vskip-6pt
\end{figure}

\noindent{\sl Pictorial representation.}
Although it is not necessary, we find it helpful to have a pictorial representation of some morphism in $\cwcob$. Since it is not easy to draw a 2-complex built on a one-point union of circles,  we will ``stretch'' its 1-skeleton and instead of  the bordism  $(X,\bar Y_{n_1},\bar Y_{n_2})$,  will draw a  CW-pair $(\tilde{X},\tilde Y_{n_1}\cup G\cup \tilde Y_{n_2})$, where 
\begin{itemize}
\item[-] $\tilde X\subset R^2\times [1,2]$ and $\tilde Y_{n_i}=\tilde X\cap R^2\times \{i\}$, $i=1,2$, is a disjoint union of $n_i$ oriented based circles;
\item[-] $G\subset \tilde{X}$ is a possibly not connected graph which intersects $R^2\times \{1,2\}$ in the base points of the components of $\tilde Y_{n_1}\sqcup\tilde Y_{n_2}$;  
\item[-] $(\tilde X,\tilde Y_{n_1}\cup G\cup \tilde Y_{n_2})/G$ is isomorphic to $(X,\bar Y_{n_1},\bar Y_{n_2})$ and the image of $G$ under the quotient map is the unique vertex of $X$, while the images of $\tilde Y_{n_1}$ and $\tilde Y_{n_2}$  are exactly  $\bar Y_{n_1}$ and $\bar Y_{n_2}$. 
\end{itemize}
 In the pictures $G$ will be drawn as a thick dark-gray {\sl unoriented} graph, while the 1-cells of $X$ will be drawn as black thick curves, which will be oriented if they belong to $\tilde Y_{n_1}$ and $ \tilde Y_{n_2}$.
 An example of a morphism $(X, \bar Y_2, \bar Y_1)$ in $\cwcob$ is presented in Figure \ref{example-cwcob/fig}. The letters $a_1, a_2, b$, etc. indicate the corresponding generators of the fundamental group of $X$ and will be used later.

\begin{figure}[htb]
  \includegraphics{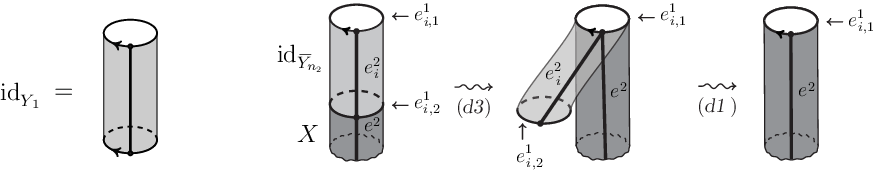}
\caption{The identity morphism.}
\label{identity-pf/fig}
%\vskip-6pt
\end{figure}

\medskip
\noindent{\sl Identity morphism.}  $\id_{\bar Y_n}$ is obtained by attaching $n$ 2-cells on $\bar Y_n\vee \bar Y_n$ as follows: we orient the boundary of   $D^2$ and identify it with the union of the two intervals lying in the upper and lower half planes, $\partial D^2=I_{1}\cup I_{2}$; then the attaching map of the $i$-th 2-cell of $\id_{\bar Y_n}$ sends  $I_{1}$ (resp. $I_2$) to the $i$-th 1-cell of the source (resp. the target) by an orientation-preserving (resp. orientation-reversing) homeomorphism ($\id_{\bar Y_1}$ is presented in  Figure \ref{identity-pf/fig}). The proof that $\id_{\bar Y_{n_2}}\circ (X,\bar Y_{n_1},\bar Y_{n_2})=(X,\bar Y_{n_1},\bar Y_{n_2})$ is illustrated in the same figure, where $e^1_{i,1}$ (resp.  $e^1_{i,2}$) is the $i$-th 1-cell in the target  (resp. source) of $\id_{\bar Y_{n_2}}$, while $e^2_i$ is its $i$-th 2-cell: we homotope in  $\id_{\bar Y_{n_2}}$ the attaching map of any 2-cell $e^2$ of $X$, which contains $e_{i,2}^1$ in its closure, until the image of such map is disjoint from the interior of  $e^1_{i,2}\cup e^2_i$ ; eventually we collapse  the pair $(e^2_{i},e^1_{i,2})$. The resulting bordism is clearly 2-equivalent to $(X,\bar Y_{n_1},\bar Y_{n_2})$.

%%%%%%%%%%%%
%%%%%%%%%%%%
\section{Group presentations and 2-dimensional CW-complexes%
\label{pres-category/sec}}

\subsection{Relative group presentations}
In Remark \ref{pres/remark} we explained how  to any standard 2-dimensional CW-complex one can associate  a  presentation of its fundamental group. In order to extend  such correspondence to standard CW-pairs, we will need the concept of a relative group presentation. In general, a relative group presentation is an expression of the form $P=\langle H; B\,|\, R\rangle$, where $H$ is a group, $B$ is a set of generators disjoint from $H$ and  $R$ is a multiset of words in $H\star \langle B\rangle$. In the present work we will use this concept in the specific case in which  $H$  is a group freely generated by a given set $A$.

\begin{definition}\label{word/def} Given a set $A$ and $n\in \N$,  a {\sl word}   of length $n$ in  $A$ is any expression of the form $w=a_1^{\epsilon_1} a_2^{\epsilon_2}\dots a_n^{\epsilon_n}$, $a_i\in A$, $\epsilon_i=\pm 1$. The exponential notation $a_i^m$, $m\in \Z-\{0\}$ will be used as an abbreviation to indicate $|m|$ consecutive terms equal to $a_i^{+1}$ if $m>0$ and to $a_i^{-1}$ if $m<0$. 

The length of the word $w$ will be indicated by $|w|$, while $|w|_{a_i}$ will indicate the number of occurrences of $a_i$ in $w$.

Observe that the concatenation defines  a product on the set of words in $A$ whose identity is the empty word and will be indicated by $1$.
\end{definition}

\begin{definition}\label{relativepr/def}
We call  {\sl relative group presentation}   an expression of the form $P=\langle A;B\;|\; R\rangle$  where $A$ and $B$ are finite sets and $R$ is a finite multiset\footnote{Two elements of $R$ represent {\sl different} relators even if they are represented by the same word. Moreover, the trivial relator $r=1$ is allowed.} of words in $G=A\sqcup B$. $A$ and $B$ will be called respectively the {\sl external} and the {\sl internal generators} of $P$, while the elements of $R$ are called the {\sl relators} of $P$. Two relative presentations $P=\langle A;B\;|\; R\rangle$ and $P'=\langle A';B'\;|\; R'\rangle$ will be considered {\sl identical} if there exist bijective maps $F_A:A\to A'$, $F_B:B\to B'$ and $F_R:R\to R'$ such that for any word $w\in R$ of length $n$, $F_R(w)$ coincides with the image of $w$ under the map $(A\sqcup B)^{\times n}\to (A'\sqcup B')^{\times n}$ induced by $F_A\sqcup F_B$. Moreover, two relative presentations will be called {\sl AC-equivalent } if one of them can be obtained from the other by a finite sequence of the following moves, to which we will refer as {\sl relative AC-moves} or shortly {\sl AC-moves}: 

\begin{itemize}
\smallskip
\item[\hrel{ac1}]\label{E:ac1} adding a new generator $b$ to $B$ and a new relator $bw$ to $R$, where $w$ is any word in $G$, or the inverse of such move;
\item[\hrel{ac2}]\label{E:ac2} substituting a relator $w\in R$ by $g g^{-1}w$ or by $g^{-1}gw$ for some  $g\in G$, or the inverse of such a move;
\item[\hrel{ac3}]\label{E:ac3} replacing a relator by its left  product with another  relator;
\item[\hrel{ac4}]\label{E:ac4} cyclic permutation of a relator;
\item[\hrel{ac5}]\label{E:ac5} replacing a relator by its inverse;
\end{itemize}
\end{definition}

\smallskip
The following proposition is well known, but we sketch the proof for completeness.

\begin{proposition}\label{add-ACmoves/thm} Two relative presentations related by each of  the following moves  are AC-equivalent:
\begin{itemize}
\smallskip
\item[\hrel{ac2'}] \label{E:ac2'} insertion or cancellation in a relator of $g g^{-1}$ or of $g^{-1} g$ for some $g\in G$;
\item[\hrel{ac3'}]\label{E:ac3'}  replacing a relator by its right  product with another  relator;
\item[\hrel{ac6}]\label{E:ac6} replacing in all relators an internal generator  $b\in B$ by $b^{-1}$, $gb$ or $bg$ for any $g\in G$, $g\neq b$ (elementary Nielsen transformation moves);
\item[\hrel{ac7}]\label{E:ac7} if  $R$ contains  a word $w'$ of the form $bw^{-1}$, $w^{-1}b$, $wb^{-1}$ or $b^{-1}w$ where $b\in B$ is an internal generator and $w$ is a word in $G\smallsetminus \{b\}$, replacing $b$ by $w$ in all relators except $w'$, and then deleting $b$ and $w'$. 
\end{itemize}
\end{proposition}

\begin{proof} \hrel{ac2'} is a generalization of \hrel{ac2} which is obtained by using \hrel{ac4} to cyclically  permute the relator in order to bring  the letter in front of which we want to insert $gg^{-1}$ to the first position, applying \hrel{ac2} and then cyclically  permuting back to the original order.

To achieve \hrel{ac3'} we  invert both relators,  replace one by its left  product with the other, and then invert back.
 
The elementary Nielsen transformation moves \hrel{ac6}  follow from \hrel{ac7} after an insertion of a new generator $c$ and a new relator respectively of the form $cb$, $cb^{-1}g$ or $cgb^{-1}$. Therefore, it remains to prove \hrel{ac7}. 

Assume now that $R=\{w^{-1}b\}\cup R'$, where $b\in B$ and $w$ is a word in $G\smallsetminus \{b\}$ (all other cases follow from this one by applying \hrel{ac4} and/or \hrel{ac5}). Then \hrel{ac7} follows by the following inductive argument on the total number of appearances of $b$ and $b^{-1}$  in the relators of $R'$. If $b$ or $b^{-1}$ don't appear in $R'$, the statement follows by \hrel{ac1}. If instead $b$ or $b^{-1}$ appear in a relator $w'\in R'$, by inverting (if necessary) and applying cyclic permutations (moves \hrel{ac4} and \hrel{ac5}) we  bring $w'$ in the form $b^{-1}z$, then we multiply it on the left by $w^{-1}b$ and cancel $bb^{-1}$ (moves \hrel{ac3} and \hrel{ac2'}), replacing in this way $w'$ by $w^{-1}z$. Eventually, we apply backwards cyclic permutations and inversions to obtain the original relator in which one of the $b$'s has been replaced by $w$. Therefore, the presentation is AC-equivalent to one in which $b$ has one fewer occurrence. 
\end{proof}

\begin{remark}\label{acmoves/rmk}
As explained in \cite[Chapter I, Section 2.3]{HMS93}, three different types of transformations of finite group presentations are relevant to the Andrews-Curtis problem. $Q$-transformations are generated by conjugating or inverting a relator and by multiplying one relator on either side by another; then $Q^*$-transformations additionally allow Nielsen transformations of the generators and finally the $Q^{**}$-transformations allow also the introduction or deletion of a generator together with the corresponding relator. 

The equivalence relation of Definition 2.2 is the relative analogue of $Q^{**}$-equivalence, with the external generators $A$ kept fixed. Indeed, as it was shown in Proposition \ref{add-ACmoves/thm}, moves \hrel{ac1}-\hrel{ac5} generate all $Q^{**}$-transformations, while the opposite is shown in \cite[Chapter I, Theorem 2.4]{HMS93}. Thus, when $A=\varnothing$, the equivalence relation of Definition 2.2 agrees with $Q^{**}$-equivalence. Throughout this paper we call it relative Andrews-Curtis equivalence, or AC-equivalence for short, following the choice of terminology made in \cite{BQ05}. We refer the interested reader to \cite{HM06} and \cite[Section 2.3]{HM17} for a discussion of the relative Andrews-Curtis problem in the more general context in which a common subcomplex is kept fixed during the transformation.

We also caution the reader that in much of the group-theoretic literature (see \cite{MMS02} and references therein), the terms Andrews-Curtis moves and Andrews-Curtis equivalence often refer to $Q$-transformations and $Q$-equivalence since the original AC-conjecture \cite{AC65} regarding contractible 2-complexes, once translated in terms of group presentations, states that any two balanced presentations of the trivial group with the same number of relators are $Q$-equivalent. 
\end{remark}

\subsection{The category $\pres$ of relative group presentations}

In what follows we will use the notation $\bar A =(a_1,a_2,\dots,a_n)$ to indicate the set $A=\{a_1,a_2,\dots,a_n\}$ whose elements have been totally ordered.

\medskip
\noindent{\sl Objects.} The  set of objects of the  category of relative presentations $\pres$ is $\N$. 

\smallskip
\noindent{\sl Morphisms.}  A morphism $\langle\, \bar A_1,\,\bar A_2;\, B\, |\, R\, \rangle: n_1\to n_2$ of the  category of relative presentations $\pres$, is given by the {\sl AC-equivalence classes} of relative presentations $\langle\,  A_1\sqcup A_2;\, B\, |\, R\, \rangle$ in which $A_i$, $i=1,2$ is a set of cardinality $n_i$, together with a choice of total order on each $A_i$. We will refer to both $n_1$ and $\bar A_1$ (resp. $n_2$ and $\bar A_2$) as the {\sl source} (resp. {\sl target}) of the above morphism. Two morphisms $\langle\, \bar A_1,\,\bar A_2;\, B\, |\, R\, \rangle, \langle\, \bar A_1',\,\bar A_2';\, B'\, |\, R'\, \rangle: n_1\to n_2$  are considered identical if the corresponding relative presentations can be identified through the bijective maps  $F_{A_1}\sqcup F_{A_2}\sqcup F_B: \bar A_1\sqcup  \bar A_2\sqcup B\to \bar  A_1'\sqcup \bar  A_2' \sqcup B'$,  where $F_{A_1}$ and $F_{A_2}$ preserve the given orders.   

\smallskip
\noindent{\sl Notation.} In the explicit expression for a given morphism in $\pres$ we will present the ordered sets $A_1$ and $A_2$ as sequences of small letters delimited by curly brackets, while the unordered sets $B$ and $R$ won't be delimited, following the usual notation for a group presentation.

\smallskip
The composition of two morphisms is defined as:
$$
\langle\, \bar A_2,\,\bar A_3;\, B'\, |\, R'\, \rangle\circ \langle\, \bar A_1,\,\bar A_2;\, B\, |\, R\, \rangle =\langle\, \bar A_1,\,\bar A_3;\, A_2\sqcup B\sqcup B'\, |\, R\sqcup R'\, \rangle,
$$
where the target of the first morphism has been identified with the source  of the second one through a bijective map which respects the sets' enumeration.
 Observe that the composition is well defined on equivalence classes since any AC-move on one of the factors induces the same move on the composition, and the associativity property is straightforward.

The identity morphism $n\to n$ is given by 
$$
\idp_n=\langle\, (a_1,a_2,\ldots, a_n)\,,(b_1,b_2,\ldots, b_n)\,;\, \emptys\, |\, a_1b_1^{-1},\,a_2b_2^{-1},\,\ldots,\, a_nb_n^{-1}\, \rangle
$$
Indeed, in order to see  that the compositions $\id_n\circ P$ and $P\circ \id_n$ (when defined) are AC-equivalent to $P$, we apply move \hrel{ac7} to substitute  each $b_i$ by $a_i$, $i=1,\dots n$, in the relators of $P$, and then delete all the relators of $\id_n$ and the generators $b_i$, $1\leq i\leq n$, from the resulting presentation.

\begin{theorem}\label{presCW/thm}
There exists a category equivalence functor $\Pi:\cwcob\to\pres$.
\end{theorem}

\begin{proof} 
On the set of objects  we define $\Pi(\bar Y_n)=n$, for any $n\in \N$. 
On the set of morphisms we define $\Pi(X,\bar Y_{n_1}, \bar Y_{n_2})=\langle\, \bar A_1,\,\bar A_2;\, B\, |\, w_1,\dots,w_s \rangle$, where $A_i$, $i=1,2$, is the set of 1-cells of $Y_{n_i}$, $B$ is the set of 1-cells of $X$ which are not in $Y_{n_1}\vee Y_{n_2}$, and  $s$ is the number of 2-cells of $X$. In order to define the $w_j$'s we fix an orientation of the elements of $B$ and for each 2-cell $e_j^2$ of $X$ with attaching map $\varphi_j:\partial D\to X^1$ choose an orientation and a base point $b_j\in\partial D$ such that $\varphi_j(b_j)=X^0$. Then
$w_j$, $j=1,\dots, s$ is the word in $A_1\sqcup A_2\sqcup B$ describing the attaching map of the $j$-th 2-cell of $X$ as an element  in $ \pi_1(X^1, X^0)$.

\smallskip
Observe  that the AC-equivalence class of $\Pi(X,\bar Y_{n_1}, \bar Y_{n_2})$ doesn't depend on the choices made in its definition. Indeed, changing the orientation of an element in $B$ replaces the corresponding generator of $\langle\, \bar A_1,\,\bar A_2;\, B\, |\, w_1,\dots,w_s \rangle$ by its inverse (move \hrel{ac6}); changing the orientation of $\partial D$ in some $\varphi_j:\partial D\to X^1$ replaces $w_j$ by its inverse (move \hrel{ac5}), while making a different choice of the base point of $\partial D$ induces a cyclic permutation of $w_j$ (move \hrel{ac4}).

\smallskip
Moreover, if $(X,\bar Y_{n_1}\vee \bar Y_{n_2})$ is changed by a 2-deformation, i.e. by some of the moves \hrel{d1-1'}, \hrel{d2-2'} or \hrel{d3}, the AC-equivalence class of $\langle\, \bar A_1,\,\bar A_2;\, B\, |\, w_1,\dots,w_s \rangle$ remains the same. Indeed, \hrel{d1-1'} adds or cancels a pair of a new generator and a new relator in the presentation as in move \hrel{ac1}; \hrel{d2-2'} inserts or cancels $g_ig_i^{-1}$ in a relator (move \hrel{ac2'}) and \hrel{d3} replaces  a relator by its product with (a cyclic permutation of) another one (moves \hrel{ac4} and \hrel{ac3}). It remains to observe that $\Pi$ preserves the composition and the identity morphisms.

\smallskip
In order to show that $\Pi$ is a category equivalence, we  construct  a functor $\pibar:\pres\to\cwcob$   such  that $\pibar\circ \Pi=\id_{\cwcob}$ and $\Pi\circ \pibar=\id_{\pres}$.
On the set of objects  $\pibar(n)=\bar Y_n$, while for any morphism $P=\langle\, \bar A_1,\,\bar A_2;\, B\, |\, R\, \rangle: n_1\to n_2$ in $\pres$, $\pibar(P)=(X_P,\bar Y_{n_1}, \bar Y_{n_2})$ where $X_P$ is the standard 2-dimensional  CW-complex built as follows. %Let $m=|B|$ and $s=|R|$ denote the cardinalities of $B$ and $R$. 
The 1-skeleton of $X_P$ is  $Y_m\vee \bar Y_{n_1}\vee \bar Y_{n_2}$, where  $m=|B|$ and the given enumeration  defines a bijective map from $A_i$ to the set of 1-cells of $Y_{n_i}$, $i=1,2$. Then by choosing an orientation of all 1-cells of $Y_m$ and a bijective map  from $B$ to the set of the 1-cells of $Y_m$,  each word $w\in R$ defines a unique element $\tilde w\in \pi_1\left(\bar  Y_m\vee \bar Y_{n_1}\vee \bar Y_{n_2}\right)$ and we complete the construction of $X_P$  by attaching  one 2-cell on each $\tilde w$. 
Observe  that different choices  leave invariant the relative 2-dimensional CW-complex $\Pi(P)$ up to cellular isomorphism. The same is true if we replace a relator of $P$ by its cyclic permutation or by its inverse (moves \hrel{ac4} and \hrel{ac5}). The rest of the AC-moves, \hrel{ac1}, \hrel{ac2} and \hrel{ac3}, change $\Pi(P)$ by the 2-deformation moves  \hrel{d1-1'}, \hrel{d2-2'} and  \hrel{d3}. Moreover, $\Pi$ respects  compositions and  identities.  

Now the identities $\pibar\circ \Pi=\id_{\cwcob}$ and $\Pi\circ \pibar=\id_{\pres}$ are straightforward.
\end{proof}

From now on we will work only with the category of relative presentations $\pres$. Nevertheless, whenever possible we will illustrate the morphisms  in $\pres$ by drawing the corresponding relative 2-dimensional CW-complexes.

%%%%%%%%%%%%%%%
%%%%%%%%%%%%%%%

\section{The category $\calg$ %
\label{Hopfc/sec}}

We recall some basic definitions and statements from the general theory of
braided monoidal categories, which are used repeatedly in the paper. The reader is referred to \cite{McL98} and \cite[Sections 2.1, 2.8, 8.1]{EGNO15} for all details.

\begin{definition}\label{strict-mon-braided-cat/def} 
A \textit{strict braided monoidal category} $\C = (\C,\diam\,,\one,\,\braid)$ is a  category $\C$ together with an associative
bifunctor $\diam: \C \times \C \to \C$, called  \textit{product} on $\C$, an object $\one$ which is right and
left unit for $\diam\,$, and a natural isomorphism $\braid_{A,B}: A \diam B \to B\diam A$, defined for any $A,B \in \Obj(\C)$, satisfying the following identities:
\begin{itemize}
\item[(i)]  $\braid_{\one,A}=\id_A=\braid_{A,\one}$;
\item[(ii)]$\braid_{A, B \diam C} = (\id_B \diam \braid_{A,C}) \circ (\braid_{A,B} \diam \id_C)$
and
$\braid_{A \diam B, C} = (\braid_{A,C} \diam \id_B) \circ (\id_A \diam \braid_{B,C})$;
\item[(iii)]  $\braid_{A,B}$ is natural in $A$ and $B$, i.e. for 
any morphisms $f: A \to A'$ and $g:B\to B'$, $\braid_{A', B}\circ(f\diam \id_B)=(\id_B\diam f )\circ\braid_{A, B}$ and $\braid_{A, B'}\circ(\id_A\diam g)=(g\diam \id_A )\circ\braid_{A, B}$.
\end{itemize}
\smallskip
A strict braided monoidal category is called {\sl symmetric} if $\braid_{B,A}\circ \braid_{A,B}=\id_{A\otimes B}$ for any $A,B \in \Obj(\C)$.
\end{definition}

All monoidal categories used in the following will be strict.

%\vspace{-15pt}
%\smallskip
%\begin{sparagraph}
\noindent{\sl Action of symmetric group.} \label{action-sym/par}
Let $\perm{n}$ be the symmetric group of degree $n\geq 1$, thought as the set of bijective functions from $\{1,2,\dots,n\}$ to itself with group operation given by the composition. For further convenience we also define $\perm{0}=\{\id_\emptyset\}$. Then, given a set $\A$, there is a standard left action of $\perm{n}$  on  $\A^{\times n}$ which sends the $i$-th element in the $\sigma(i)$-th place:
$$
\sigma(A_1,A_2,\dots ,A_n)=(A_{\sigma^{-1}(1)} ,A_{\sigma^{-1}(2)}, \dots, A_{\sigma^{-1}(n)}).
$$
 We will use this action in different contexts. The first one is when $\A=\Obj(\C)$ for some monoidal category $(\C,\diam)$ and $\A^{\times n}$ is identified with the set of $\diam$-words $A_1\diam A_2\diam\dots\diam A_n$ in $\A$. The second one is the action of $\perm{n}$ on the set of words of length $n$ in some set $G$ when  $\A=G\times \{+1,-1\}$ (see Definition \ref{word/def}).
 %\end{sparagraph}
 
 \begin{definition}\label{iota/defn}
Given any $n\geq 1$ and any decomposition $n=n_1+n_2+\dots n_k$,  $\iota_{n_1,\dots,n_k}:\times_{i=1}^k \perm{n_i}\hookrightarrow \perm{n}$ is the injective group homomorphism which applies $\sigma_i\in\perm{n_i}$ on the $n_i$-subtuple $(n_1+\dots +n_{i-1}+1,\dots, n_1+\dots +n_{i-1}+n_i)$ of the $n$-tuple $(1,2,\dots, n)$.
\end{definition}
In what follows, whenever the context is clear, we will identify $\times_{i=1}^k \perm{n_i}$ with the image of $\iota_{n_1,\dots,n_k}$ writing $ \sigma_1\times \sigma_2\times\dots\times \sigma_k$ instead of $ \iota_{n_1,\dots,n_k}(\sigma_1,\sigma_2,\dots,\sigma_k)$.

\medskip
We  now  state the coherence theorem of Mac Lane for symmetric categories, specified to the case of  {\sl strict} monoidal categories, which basically says that any  morphism from $A_1\diam A_2\diam\dots\diam A_n$ to $\sigma(A_1\diam A_2\diam\dots\diam A_n)$, built up as a composition of products of $\braid$'s and identities, is uniquely determined by the permutation $\sigma\in \perm{n}$. 

\begin{theorem}[Theorem 1, Chapter XI, \cite{McL98}] \label{symm-coherence/thm}
Let $(\C,\diam, \one, \braid)$ be a symmetric strict monoidal category and $A_1, \dots, A_n\in \Obj(\C)$. For any permutation $\sigma\in \perm{n}$ there exists a unique 
 morphism $\Upsilon^{\sigma}_{A_1, \dots, A_n}: A_1\diam A_2\diam\dots\diam A_n\to \sigma(A_1\diam A_2\diam\dots\diam A_n)$, natural in the $A_i$'s, such that
 
\begin{itemize}
\item[(i)] $\Upsilon^{\id}_{A_1, \dots, A_n}=\id_{A_1\diam A_2\diam\dots\diam A_n}$;
\item[(ii)] $\Upsilon^{(1,2)}_{A_1, A_2}=\braid_{A_1, A_2}$, where $(1,2)\in \perm{2}$ denotes the transposition;
\item[(iii)] $\Upsilon^{\sigma'}_{A_{\sigma^{-1}(1)},\dots,A_{\sigma^{-1}(n)}}\circ\Upsilon^{\sigma}_{A_1, \dots, A_n}=\Upsilon^{\sigma'\circ\sigma}_{A_1, \dots, A_n}$;
\item[(iv)] $\Upsilon^{\sigma}_{A_1, \dots, A_n}\diam\Upsilon^{\sigma'}_{A_1', \dots, A_m'}=\Upsilon^{\sigma\times\sigma'}_{A_1, \dots, A_n,A_1', \dots, A_m'}$ for every $\sigma\in \perm{n},\;\sigma'\in\perm{m}$;
%\item[(v)] $\Upsilon^{\sigma}_{A_1, \dots, A_n}$ is natural in the $A_i$'s, i.e. given any morphism $F:A_i\to B$, 
%$$(\id_{\sigma(i)-1}\diam F\diam\id_{n-\sigma(i)})\circ \Upsilon^{\sigma}_{A_1, \dots, A_n}=\Upsilon^{\sigma}_{A_1, \dots,A_{i-1} ,B, A_{i+1},\dots, A_n}\circ (\id_{i-1}\diam F\diam\id_{n-i})$$
\end{itemize}

$\Upsilon^{\sigma}_{A_1, \dots, A_n}$ is called {\sl canonical morphism}. 
\end{theorem}

\begin{definition}\label{autonomous/def}
A  strict braided monoidal category $\C = (\C,\diam\,,\one,\braid)$ is called
\textit{rigid} (or autonomous) %(\cite[Sections 2.10]{EGNO15}, \cite{Sh94}) 
if for every $A \in \Obj \C$ it is given
a \textit{right dual} object $A^\ast \in \Obj \C$ and two morphisms
$$
\begin{array}{ll} 
  \ev_A: A \diam A^\ast \to \one \ \ &\text{(evaluation)\,,}\\[2pt]
\coev_A: \one \to A^\ast \diam A \ \ & \text{(coevaluation)\,,}
\end{array}
$$
such that
$ (\ev_A \diam \id_A)\circ (\id_A \diam \coev_A)=\id_A\quad\text{and}\quad 
(\id_{A^\ast} \diam \ev_A)\circ (\coev_A \diam \id_{A^\ast})=\id_{A^\ast}$.

Then, given any morphism $F: A \to B$ in $\C$, we define its 
\textit{dual} $F^\ast: B^\ast \to A^\ast$ as follows:
$$F^\ast = (\id_{A^\ast} \diam \ev_B) \circ (\id_{A^\ast} \diam F \diam \id_{B^\ast})
\circ (\coev_A \diam \id_{B^\ast})\,.$$
\end{definition}

\begin{definition}\label{hopf-algebra/def}
Given a  strict braided monoidal category $(\C, \diam,\one, \braid)$, a \textit{braided bialgebra}
in $\C$ is an object $H$  equipped with the following
families of morphisms  (cf. Table \ref{axiomsH/fig}):
\begin{itemize}
 \item a \textit{product} $\prodH : H \otimes H \to H$ and a \textit{unit} $\unitH : \one \to H$;
 \item a \textit{coproduct} $\coprH : H \to H \otimes H$ and a \textit{counit} $\counH : H \to \one$,
\end{itemize}
subject to the following axioms:
\begin{gather*}
\m \circ (\m \diam \id_H) = \m \circ (\id_H \diam \m),
 \tag*{\hrel{a1}}
 \\
\m\circ(\eta\diam\id_H)  = \id_H = \m\circ(\id_H\diam\eta),
 \tag*{\hrel{a2-2'}}
 \\
 (\Delta\diam\id_H) \circ \Delta = (\id_H\diam\Delta) \circ \Delta;
 \tag*{\hrel{a3}}
 \\
 (\epsilon\diam\id_H)\circ\Delta = \id_H = (\id_H\diam\epsilon)\circ\Delta;
 \tag*{\hrel{a4-4'}}
 \\
 (\m \diam \m) \circ (\id_H \diam \braid_{H,H}\diam \id_H) \circ
(\Delta \diam \Delta) = \Delta \circ \m,
 \tag*{\hrel{a5}}\\
 \epsilon \circ \m = \epsilon \diam \epsilon;                 \tag*{\hrel{a6}}
\\
 \Delta \circ \eta = \eta \diam \eta,
 \tag*{\hrel{a7}}\\
 \epsilon \circ \eta = \id_{\one};
 \tag*{\hrel{a8}}
 \end{gather*}

%\noindent 
A bialgebra  $H$  in $\C$ is  
\begin{itemize}
\item[-] \textit{cocommutative } if
\begin{gather*}\braid_{H,H}\circ \Delta=\Delta,
 \tag*{\hrel{a9}}\end{gather*}

\item[-]
 \textit{unimodular}
 \footnote{In the literature the term \textit{unimodular} is not uniquely defined. In \cite{Mo93} a Hopf algebra is called unimodular if it has a two-sided integral element, while in \cite{Mj95} it is called unimodular if it has a two-sided integral form. Our definition is stronger since it implies both and adds the normalization condition \hrel{i3}.} 
 if it is equipped with a
 \textit{two-sided integral form } $\lambda : H \to \one$  and a \textit{two-sided integral element }\footnote{Here we follow the terminology in \cite{BBDP23}, but in  the usual Hopf-algebra terminology (see for example \cite{Mo93}), $\Lambda:\one\to H$ is called integral {\sl in} $H$ and $\lambda:H\to\one$ an integral {\sl on} $H$. On the other hand, in the Quantum Topology literature, $\Lambda$ is usually called {\sl cointegral} of $H$, while $\lambda$ is called {\sl integral} of $H$. }  $\Lambda :  \one\to H$ such that
\vskip-12pt
\begin{gather*}(\id_H \diam \lambda) \circ \Delta = \eta \circ \lambda= (\lambda \diam \id_H) \circ \Delta,
 \tag*{\hrel{i1-1'}}\\
 \m \circ (\Lambda \diam \id_H) = \Lambda \circ \epsilon= \mu \circ (\id_H \diam \Lambda) ,
 \tag*{\hrel{i2-2'}}\\
 \lambda \circ \Lambda = \id_\one,
 \tag*{\hrel{i3}}
 \end{gather*}

\item[-]
\textit{Hopf}  if it is equipped with an  invertible morphism $S : H \to H$, called the \textit{antipode}  of $H$, such that
\vskip-5pt
\begin{gather*}\m\circ(S\diam\id_H)\circ\Delta = \eta\circ\epsilon=\m\circ(\id_H\diam S)\circ\Delta.
 \tag*{\hrel[s1]{s1-1'}}\\
 \antipH \circ \antipH^{-1} = \id = \antipH^{-1} \circ \antipH.
 \tag*{\hrel[s2]{s2-3}}\end{gather*}
\end{itemize}
\end{definition}

\begin{table}[p]
  \includegraphics{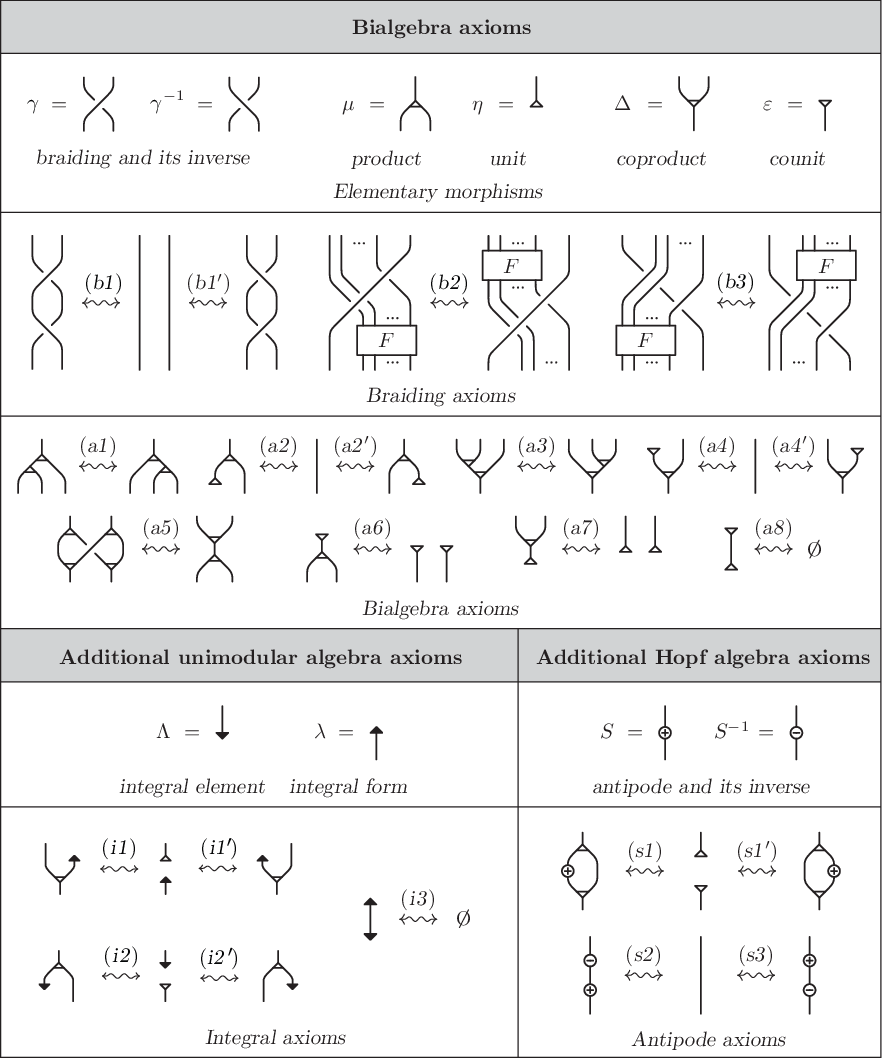}
\caption{}
\label{axiomsH/fig}
 \label{E:i1}\label{E:i1'} \label{E:i1-1'} \label{E:i2} \label{E:i3} \label{E:i2'} \label{E:i2-2'} \label{E:a3} \label{E:a1} \label{E:a2'}\label{E:a2-2'} \label{E:a2}  \label{E:a4-4'}\label{E:a4} \label{E:a4'}\label{E:a5} \label{E:a6} \label{E:a7} \label{E:a8}
 \label{E:s1} \label{E:s1'} \label{E:s2}\label{E:s3} 
 \label{E:b1} \label{E:b1'} \label{E:b2} \label{E:b3} 
%\vskip-3pt
\end{table}

\begin{table}[p]
  \includegraphics{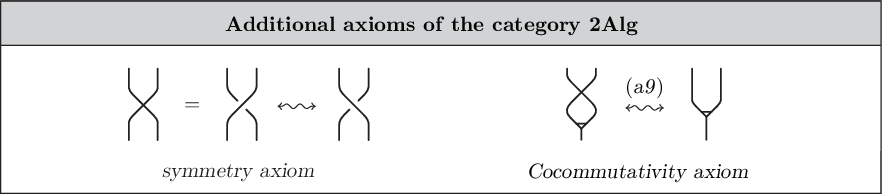}
\caption{}
\label{axioms2Alg/fig}
 \label{E:a9} 
%\vskip-3pt
\end{table}

\begin{definition}\label{calg/def}
Let $\calg$ be the strict symmetric monoidal category freely generated by a unimodular cocommutative   bialgebra $H$ (see Tables \ref{axiomsH/fig} and \ref{axioms2Alg/fig}). In particular,   if $\C$ is any symmetric monoidal category with a unimodular cocommutative   bialgebra $H_\C$ in it, then there exists a unique strict symmetric monoidal functor $\calg \to \C$ sending $H$ to $H_\C$.
\end{definition}

\begin{table}[t]
  \includegraphics{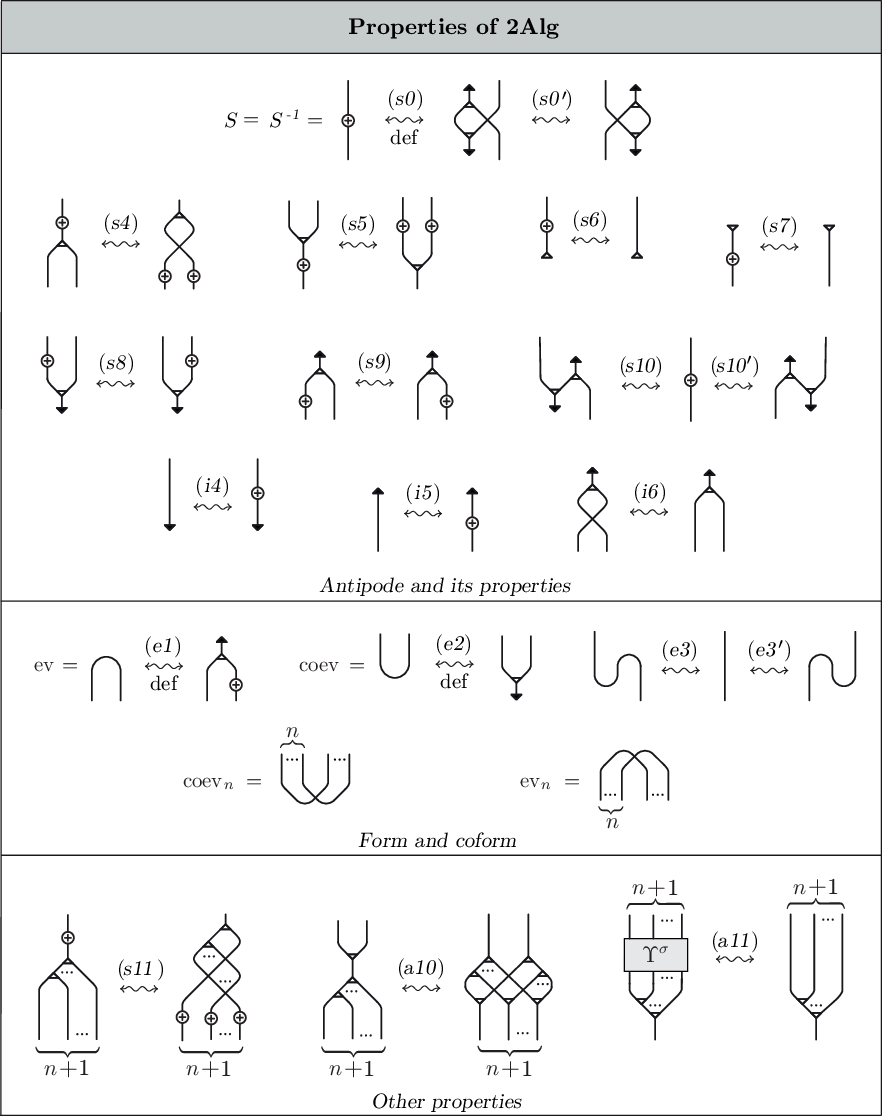}
\caption{}
\label{antipode-form/fig}
\label{E:s0}\label{E:s0'} \label{E:s10} \label{E:s10'}  \label{E:s5} \label{E:s4} \label{E:s4-5} \label{E:s6}  \label{E:s7} \label{E:s6-7} \label{E:s8} \label{E:s9} \label{E:s11} \label{E:i4} \label{E:i5} \label{E:i4-5} \label{E:i6}\label{E:e1} \label{E:e2} \label{E:e3} \label{E:a10} \label{E:a11}
\end{table}

%\begin{sparagraph}
\noindent{\sl Diagrammatic language.} The set of objects of $\calg$ is the free monoid generated by $H$, while the set of morphisms consists of all the compositions of products of identities and of the elementary morphisms $\braid,  \Delta,\epsilon, \m, \eta,   \Lambda, \lambda$ modulo the defining axioms for a braided symmetric structure and  the defining axioms listed in Definition \ref{hopf-algebra/def}. Therefore, applying the Penrose graphical notation (see \cite{Tu94, TV17}), we describe  $\calg$  as the category of planar diagrams in $[0,1] \times [0,1]$, where the objects are sequences of points in $[0,1]$  and the morphisms are iterated products and compositions of the elementary diagrams presented in Table \ref{axiomsH/fig}, modulo the relations presented there and also plane isotopies which preserve the $y$-coordinate. Recall that the composition of diagrams $D_2 \circ D_1$ is obtained by stacking $D_2$ on the top of $D_1$ and then rescaling, while the product $D_1 \diam D_2$ is given by the horizontal juxtaposition of $D_1$ and $D_2$ and rescaling. 

%\smallskip
Most of the proofs here will consist of showing that some morphisms in the universal algebra are equivalent, meaning that the graph diagram of one of them can be obtained from the graph diagram of the other by applying a sequence of the defining relations (moves) of the algebra axioms. We will outline the main steps in this procedure by drawing in sequence some intermediate diagrams, and for each step we will indicate in the corresponding order, the main moves needed to transform the diagram on the left into the one on the right. Notice that the moves represent equivalences of diagrams and we use the same notation for them and their inverses. 
%\end{sparagraph}

\smallskip
From now on $H$ will denote the unimodular cocommutative   bialgebra which generates $\calg$. 

\medskip
%\begin{sparagraph}
The notations  \label{notations/par}  for the identities, the braiding and the canonical morphisms in Theorem \ref{symm-coherence/thm} on multiple products of $H$,  will be simplified as follows ($n,m\geq 0$, $\sigma \in \perm{n}$):  
\begin{gather*}
H^0\defeq\one, \quad H^n\defeq H^{\diam n}, \;n\geq 1\\
 \id_n\defeq \id_{H^{n}}, \qquad \braid_{n,m}\defeq\braid_{H^{n},H^{m}},\qquad \Upsilon^{\id_\emptyset}\defeq \id_\one,\quad\Upsilon^{\sigma}\defeq\Upsilon^\sigma_{H,\dots,H}:H^{n}\to H^{n}, \;n,m\geq 1.
 \end{gather*}
%Moreover, in the following we will repeatedly use the canonical morphisms associated with the following  specific families of permutations:
%\begin{itemize}
%\item[-] the permutation $\pi_{n}\in\perm{2n}$, which  brings all odd indices in front of the even ones, preserving the ordering in each of the two subsets:
%$$
%\pi_{0}=\id_\emptyset, \qquad\pi_{n}(2i-1)=i, \; \pi_{n}(2i)=n+i,\;  1\leq i\leq n;
%$$  
%\item[-]  the permutation $\tau_{n}\in\perm{n}$, which reverses the order of the indices in the $n$-tuple: 
%$$
%\tau_{0}=\id_\emptyset,\qquad \tau_{n}(i)=n+1-i,\;1\leq i\leq n.
%$$
%\end{itemize}
%\medskip
The next proposition introduces some  properties of $\calg$, the most important of which is the existence of an antipode. 

\begin{proposition} \label{antipode/thm}  Let $H$ be a unimodular cocommutative  bialgebra in $(\C, \diam,\one, \braid)$. Then (see Table \ref{antipode-form/fig})
 $H$ is a Hopf algebra with antipode 
\begin{gather*}S=((\lambda\circ\m)\diam\id_1 )\circ(\id_1\diam\braid)\circ( (\Delta\circ \Lambda)\diam\id_1)=(\id_1\diam (\lambda\circ\m))\circ(\braid\diam\id_1)\circ(\id_1 \diam(\Delta\circ \Lambda)) \tag*{\hrel[s0]{s0-0'}}\end{gather*}
Moreover, $S^2=\id_1$ and
\begin{gather*} %S^2=\id_1 \tag*{\hrel{s2}}\\
S \circ \m = \m \circ (S \diam S) \circ \braid_{1,1}\,,
 \tag*{\hrel{s4}}\\
\Delta \circ S = (S \diam S)  \circ\braid_{1,1}\circ \Delta= (S \diam S)\circ \Delta\,,
 \tag*{\hrel{s5}}\\
  S \circ \eta = \eta \,, \qquad \epsilon \circ S = \epsilon,
 \tag*{\hrel{s6-7}}\\
 (S\diam \id)\circ\Delta\circ \Lambda = (\id\diam S)\circ\Delta\circ \Lambda \,,
 \tag*{\hrel{s8}}\\
 \lambda\circ \m \circ (S\diam \id)= \lambda\circ \m \circ (\id\diam S) \,,
 \tag*{\hrel{s9}}\\
 S\circ \Lambda=\Lambda,\qquad \lambda\circ S=\lambda  \tag*{\hrel{i4-5}}\\
 \lambda\circ \m \circ \braid= \lambda\circ \m  \,,
 \tag*{\hrel{i6}}\\
 (\id\diam (\lambda\circ\m))\circ ((\Delta\circ \Lambda)\diam\id)=S=((\lambda\circ\m)\diam\id )\circ (\id\diam (\Delta\circ \Lambda)).
 \tag*{\hrel[s10]{s10-10'}}\end{gather*}
\end{proposition}

\begin{figure}[htb]
  \includegraphics{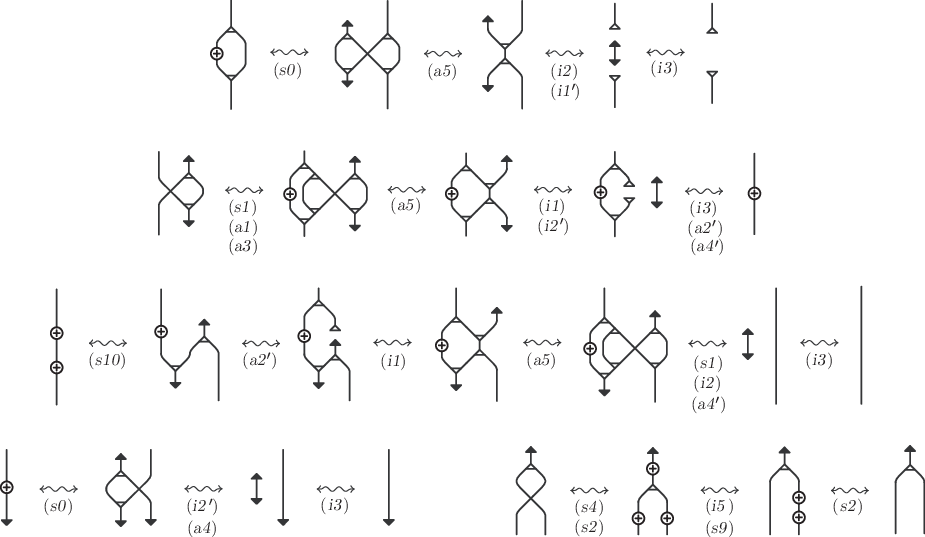}
\caption{Proof of Proposition \ref{antipode/thm}}
\label{proof-antipode/fig}
%\vskip-6pt
\end{figure}

\begin{proof} 
The proof is presented  in  Figure \ref{proof-antipode/fig}. We start by proving in the first line in the figure that $S$ satisfies  the defining property  \hrel{s1} of the antipode. Then in the second line we show  that the antipode can also be expressed by the formula \hrel{s0'}, symmetric to \hrel{s0}, and this implies \hrel{s1'} by symmetry. Now, by applying the cocommutativity condition \hrel{a9} ($\Delta=\braid\circ\Delta$) to the comultiplication morphisms in \hrel[s0]{s0-0'}, we obtain \hrel[s10]{s10-10'} and, expressing the antipode through \hrel[s10]{s10-10'}, we see that the right- and the left-hand side of \hrel{s8} and \hrel{s9} are the same. 

The fact that $S^2=\id_1$ is proven in the third line in Figure \ref{proof-antipode/fig} and \hrel{i4} is proven in the fourth. Then the proof of \hrel{i5} is obtained by reflecting the proof of \hrel{i4} with respect to the horizontal line. 

Finally, \hrel{s4}--\hrel{s7} are standard properties of the antipode (see \cite{Mo93} or \cite{BBDP23} for the diagrammatic proof), while \hrel{i6} follows from \hrel{s4} and \hrel{i5} as it is shown in the last line in Figure \ref{proof-antipode/fig}.
\end{proof}

Proposition \ref{antipode/thm} implies that the category $\calg$ (Definition \ref{calg/def}) can be equivalently defined as the symmetric monoidal category freely generated by a unimodular cocommutative Hopf algebra.

%\newpage

%%%
\begin{remark}\label{ant/rmk}
Proposition \ref{antipode/thm} appears in different contexts in the literature. First of all, it is closely related to the Larson-Sweedler theorem, according to which a finite-dimensional bialgebra with a nonsingular integral is a Hopf algebra \cite{LS69}. A categorical version of this theorem is presented in \cite{BFVV21}. The theory of integrals for Hopf algebras in rigid braided monoidal categories is developed in \cite{BKLT00}, where, however, the antipode is part of the initial Hopf-algebra structure. Finally, formulas \hrel[s0]{s0-0'} and \hrel[s10]{s10-10'} are the cocommutative unimodular specializations of Kuperberg's formulas expressing  the antipode and its inverse in terms of the integrals \cite[Lemma 3.4 and the discussion following Lemma 3.7]{Ku96}.
%The assumptions here differ from those in the usual Larson-Sweedler theorem. We do not assume that $H$ is finite-dimensional, linear, or dualizable. Instead, we assume the existence of a specified two-sided integral and cointegral satisfying$$ \lambda\circ\Lambda=\operatorname{id}_{\one}. $$ Proposition \ref{lambda/thm} then shows that this datum makes $H$ self-dual. 
\end{remark}

\begin{proposition}\label{lambda/thm}
$\calg$ is a rigid category with  $(H^n)^\ast =H^n$ for any $n\geq 0$ and (see Table \ref{antipode-form/fig})
\begin{gather*}
\ev_0=\coev_0=\id_\one\\
\ev_n=(\lambda\circ\mu\circ(\id\diam S))^{\diam n}\circ\Upsilon^{\pi^{-1}_n}: H^{n}\diam H^{n}\to \one,
\\
\coev_n=\Upsilon^{\pi_n}\circ \Delta^{\diam n}\circ \Lambda^{\diam n}:\one\to H^{n}\diam H^{n},\;n\geq 1,
\end{gather*}
where $\pi_{n}\in\perm{2n}$ is the permutation which  brings all odd indices in front of the even ones as follows:
$$
\pi_{0}=\id_\emptyset, \qquad\pi_{n}(2i-1)=i, \; \pi_{n}(2i)=n+i,\;  1\leq i\leq n.
$$  
\end{proposition}

\begin{proof}
Relations \hrel[s8]{s8-9} and \hrel[s10]{s10-10'} in Proposition \ref{antipode/thm} imply that  (see Table \ref{antipode-form/fig}):
\vskip-12pt
\begin{gather*}
 \ev =\ev_1= \lambda \circ \m \circ (\id_1 \diam S)
\,, \tag*{\hrel{e1}}\\
\coev=\coev_1  = \Delta \circ \Lambda \,,
 \tag*{\hrel{e2}}
\end{gather*}
satisfy the identities:
\begin{gather*}
(\id_1\diam \ev)\circ (\coev\diam\id_1)=\id_1=(\ev\diam\id_1)\circ(\id_1\diam \coev) \tag*{\hrel[e3]{e3-3'}}\end{gather*}
For $n> 1$, $\ev_n: H^{n}\diam H^{n}\to \one$ and  $\coev_n:\one\to H^{n}\diam H^{n}$  can be obtained by the following inductive definition\footnote{ For convenience the choice  of the  evaluation and  coevaluation  morphisms here  is different from the one  in \cite{BP11, BBDP23}.}:
\begin{gather*}
\coev_n= (\id_{n-1}\diam\braid_{n-1,1}\diam \id_1)\circ (\coev_{n-1}\diam\coev)\,,
\\
\ev_n=  (\ev_{n-1}\diam\ev)\circ(\id_{n-1}\diam\braid_{1,{n-1}}\diam \id_1).
\end{gather*}
Then the naturality of the braiding morphisms and \hrel[e3]{e3-3'} imply that  the zig-zag condition holds for any $n\geq 1$.
\end{proof}

\begin{definition}\label{deltan/defn}
For any $n\geq -1$, let $\Delta^n: H\to H^{(n+1)}$ and $\m^n:H^{(n+1)}\to H$ be  defined inductively as follows:
\begin{gather*}
 \Delta^{-1}=\epsilon, \;\;\Delta^0=\id_1 \;\text{ and }\; \Delta^n=(\Delta^{n-1}\diam\id_1)\circ\Delta\quad \text{for $n>0$},\\
\m^{-1}=\eta, \;\;\m^0=\id_1\;\text{ and }\; \m^n=\m\circ (\m^{n-1}\diam\id_1)\quad \text{for $n>0$.}
\end{gather*}
\end{definition}

The next proposition generalizes to any $n$ relations \hrel{a5}, \hrel{a9} and \hrel{s4} in Tables \ref{axiomsH/fig} and \ref{antipode-form/fig}.

\begin{proposition}\label{delta-mu-n/thm}
For any $n\geq 0$ (see Table \ref{antipode-form/fig}):
\begin{gather*} 
\Delta\circ\mu^n=(\mu^n\diam\mu^n)\circ \Upsilon^{\pi_{n+1}}\circ \Delta^{\diam (n+1)}.  \tag*{\hrel{a10}}\\
\Upsilon^\sigma\circ \Delta^{n}=\Delta^{n}, \;\sigma \in \perm{n+1} \tag*{\hrel{a11}}\\
S\circ\mu^n=\mu^n\circ \Upsilon^{\tau_{n+1}}\circ S^{\diam (n+1)},  \tag*{\hrel{s11}}
\end{gather*} 
where $\tau_{n}\in\perm{n}$ is the permutation which reverses the order of the indices in the $n$-tuple as follows: 
$$
\tau_{0}=\id_\emptyset,\qquad \tau_{n}(i)=n+1-i,\;1\leq i\leq n.
$$
\end{proposition}

\begin{proof}
\hrel{a10} is proved  by induction on $n$. If $n=0$, the statement is trivial and for $n=1$, it coincides with \hrel{a5}. Suppose now the statement is true for some $n\geq 1$. Then
\beqn
&\Delta\circ\mu\circ(\mu^n\diam\id_1)&=
(\mu\diam\mu)\circ(\id_1\diam \braid_{1,1}\diam\id_1)\circ((\Delta\circ\mu^n)\diam\Delta)\\
&&=(\mu\diam\mu)\circ(\id_1\diam \braid_{1,1}\diam\id_1)\circ(((\mu^n\diam\mu^n)\circ\Upsilon^{\pi_{n+1}})\diam\id_2)\circ\Delta^{\diam(n+2)}\\
&&=(\mu^{n+1}\diam\mu^{n+1})\circ (\id_{n+1}\diam \braid_{n+1,1}\diam\id_1)\circ(\Upsilon^{\pi_{n+1}}\diam\id_2)\circ\Delta^{\diam(n+2)}\\
&&=(\mu^{n+1}\diam\mu^{n+1})\circ \Upsilon^{\pi_{n+2}}\circ \Delta^{\diam (n+2)},
\eeqn
where we have used first  \hrel{a5}, then the induction hypothesis, and finally the naturality of the braiding.

\smallskip
For $n=0$ \hrel{a11} follows from the definitions of $\Upsilon^{\id_1}$ and that of $\Delta^0$. On the other hand, if $n>0$, any permutation in $\perm{n+1}$ can be expressed as a composition of adjacent transpositions and, according to property (iii) in Theorem \ref{symm-coherence/thm}, $\Upsilon^{\sigma'\circ\sigma}=\Upsilon^{\sigma'}\circ\Upsilon^{\sigma}$. Therefore, it is enough to prove \hrel{a11} in the case in which $\sigma$ is an adjacent transposition and in that case it follows directly from the coassociativity axiom \hrel{a3} and the cocommutativity axiom \hrel{a9}.

\smallskip
\hrel{s11} follows by induction on $n$. For $n=0$ the statement is trivial. If it is true for $n\geq 0$, then using the property of the antipode \hrel{s4} and the naturality of the braiding morphism \hrel{b3}, we obtain:
\beqn
&S\circ(\mu\circ(\mu^n\diam\id_1))&=\mu\circ\braid_{1,1}\circ((S\circ\mu^n)\diam S)\\
&&=\mu\circ\braid_{1,1}\circ((\mu^n\circ \Upsilon^{\tau_{n+1}}\circ S^{\diam (n+1)})\diam S)\\
&&=\mu^{n+1}\circ \Upsilon^{\tau_{n+2}}\circ S^{\diam (n+2)}.
\eeqn
\end{proof}

\begin{remark}
    Any cocommutative  unimodular Hopf algebra in a symmetric category is a commutative Frobenius algebra with Frobenius coproduct $\Delta_F=\Delta$, counit $\epsilon_F=\epsilon$ and Frobenius product and unit are given by 
    $$\mu_F=(\id_1\otimes \ev)\circ (\Delta\otimes \id_1)=( \ev\otimes \id_1)\circ ( \id_1\otimes \Delta),\;\; \eta_F=\Lambda.$$ The associativity and the Frobenius   relation are  consequences of the identity between the two expressions for $\mu_F$ given above and the proof of this identity in the more general context of a BP Hopf algebra can be found in \cite[Section 3.4]{BBDP23}. On the other hand, the fact that $\eta_F=\Lambda$ follows directly from the zig-zag identities \hrel[e3]{e3-3'}. 
    
This construction is a categorical version of the classical Frobenius structure associated with a finite-dimensional Hopf algebra and its integrals \cite{LS69, Pa71}. For analogous constructions from a normalized pair consisting of an integral form and an integral element in symmetric or rigid monoidal categories, see \cite{CKQW25}.
\end{remark}

%%%%%%%%%%%%%%%%%%%%%%%%%%%

\section{The equivalence functor $\Psi:\calg\to \pres$}
\label{hopf-pres/sec}

This section is dedicated to the proof of Theorem \ref{main1/thm}.

\smallskip
\begin{sparagraph}{Monoidal structure of $\pres$} \label{monoid-pres/par}  The category of relative group presentations $\pres$ carries a strict monoidal structure, whose product $\diam:\pres\times \pres\to\pres$ is given by:
\beqn
&& n\diam m=n+m \\
&&\langle\, \bar A_1,\,\bar A_2;\, B\, |\, R\, \rangle\diam \langle\, \bar A_1',\,\bar A_2';\, B'\, |\, R'\, \rangle=\langle\, \overline{ A_1\sqcup  A_1'},\,\overline{A_2\sqcup  A_2'};\,  B\sqcup B'\, |\, R\sqcup R'\, \rangle,
\eeqn
where in $\overline{ A_i\sqcup  A_i'}$,  the enumeration of $A_i'$ is shifted forward by the cardinality of $A_i$, $i=1,2$. The unit object and morphism are respectively 0 and the empty presentation. Observe that the product is well defined on AC equivalence classes and under the equivalence functor in Theorem \ref{presCW/thm} it corresponds to one-point union of bordisms in $\cwcob$.
\end{sparagraph}

\begin{table}[p]
  \includegraphics{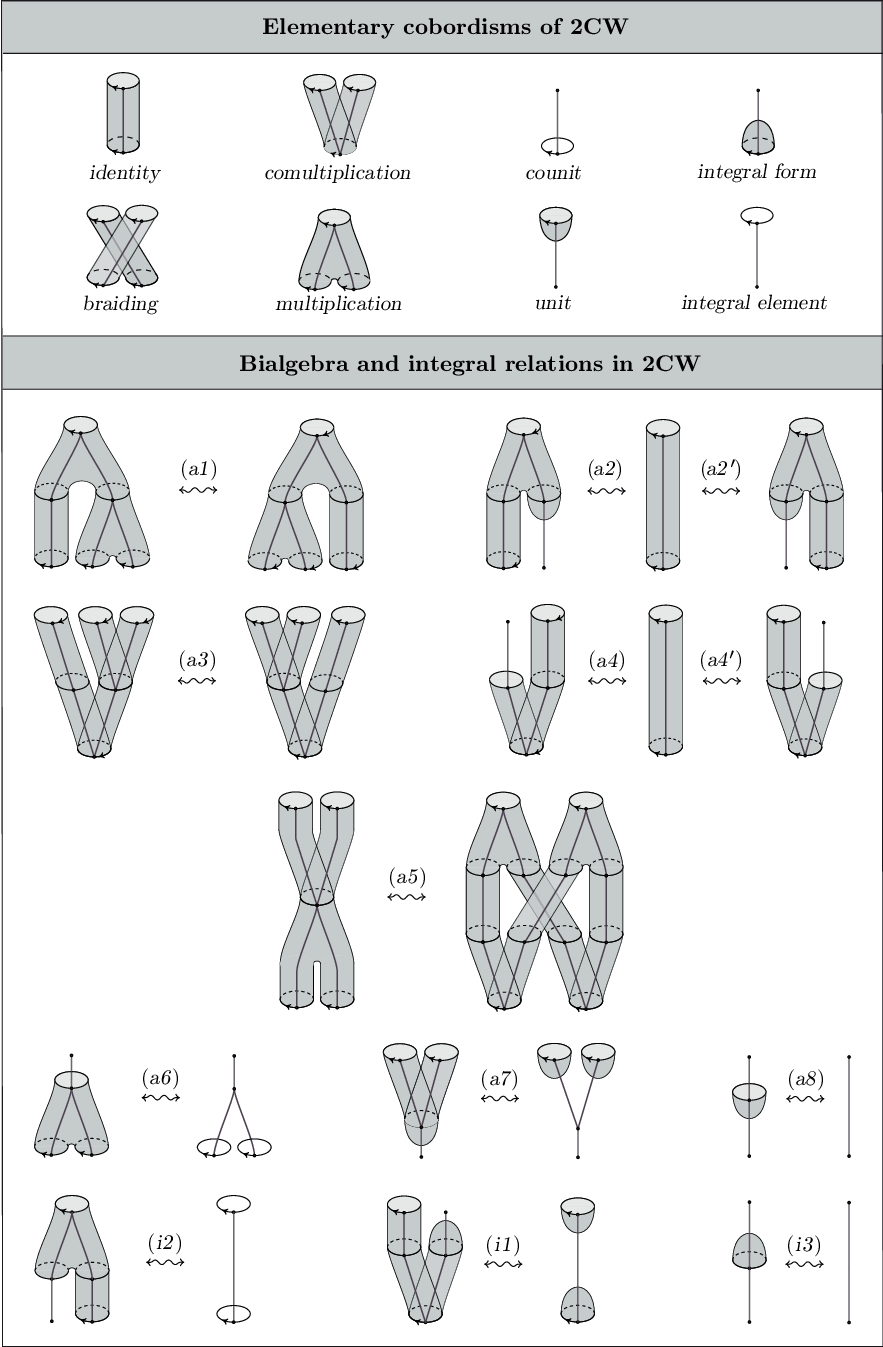}
\caption{}
\label{table-bialgebra-cw/fig}
%\vskip-3pt
\end{table}

%\medskip
\begin{theorem}\label{Hpres/thm}
$(\pres, \diam)$ is a strict symmetric monoidal category with  braiding morphisms $\gammap_{n,m}:n\diam m\to m\diam n$ given by
$$
\gammap_{n,m}=\langle\, (a_1,\dots,a_{n+m}),(b_1,\dots,b_{n+m});\emptys\, |\, a_1b_{m+1}^{-1}\, \dots, a_nb_{m+n}^{-1},a_{n+1}b_1^{-1},\dots,a_{m+n}b_m^{-1}\rangle.
$$ 
Moreover, the object $1\in\N$ is a unimodular cocommutative  Hopf algebra in $\pres$ with 
\begin{itemize}
%\item[] identity $\idp_n=\langle(a_1,\dots,a_n),(b_1,\dots,b_n);\emptys |\,a_1b_1^{-1}, \dots, a_nb_n^{-1}\rangle, : n \to n$;
%\smallskip
\item[] comultiplication $\Deltap=\langle (a),(b,c);\emptys\,|\,ab^{-1},ac^{-1}\rangle\, : 1 \to 2$;
\smallskip
\item[] counit $\epsp=\langle (a),\emptys;\emptys\,|\,\emptys\rangle\, : 1 \to 0$; 
\smallskip
\item[] multiplication $\multp=\langle (a,b),(c);\emptys\,|\,abc^{-1}\rangle\, : 2 \to 1$;
\smallskip
\item[] unit $\etap=\langle \emptys, (a);\emptys\,|\,a\rangle\, : 0 \to 1$; 
\smallskip
\item[] integral form  $\lp=\langle (a), \emptys;\emptys\,|\,a\rangle\, : 1 \to 0$; 
\smallskip
\item[] integral element  $\Lp=\langle \emptys,(a);\emptys\,|\,\emptys\rangle\, : 0 \to 1$;
\item[] antipode $\Sp=\langle (a), (b);\emptys\,|\,ab\rangle\, : 1 \to 1$.
%\smallskip 
\end{itemize}
\end{theorem}
%\vskip-15pt
The statement of Theorem \ref{Hpres/thm} is illustrated in Table  \ref{table-bialgebra-cw/fig}, where we have represented the relative 2-dimensional CW-complexes, whose images under the equivalence functor $\Pi$ in Theorem \ref{presCW/thm} are $\idp_1$, $\Deltap$, $\multp$, $\epsp$, $\etap$, $\Lp$ and $\lp$ and   the relations in $\cwcob$ which correspond to the bialgebra and integral axioms.

\begin{proof}
We have to check that the morphisms in the statement of the theorem satisfy the axioms in Definition \ref{hopf-algebra/def} (cf. Table \ref{axiomsH/fig}). All proofs are based on the following observation: in the composition $P_2\circ P_1$ of two morphisms in $\pres$, the external generators in the target of $P_1$ (and in the source  of $P_2$) become internal, and those of them which appear in a relator $r$ once with exponent $\pm 1$, can be eliminated by first cyclically  permuting and/or inverting $r$ (moves \hrel{ac4} and \hrel{ac5}) in order  to represent it as a product of some word $w$ and the inverse of the generator, then substituting the generator by $w$ in all other relations and finally canceling the generator and $r$ from $P_2\circ P_1$ (move \hrel{ac7}). Using this strategy we will  eliminate all internal generators which appear in both the right- and left-hand sides of the relations in  Definition \ref{hopf-algebra/def} and show that the resulting presentations are equivalent. 
 
%\medskip
We start by showing that the braiding is symmetric. Indeed,
\beqn
&\gammap_{1,1}\circ \gammap_{1,1}&=\langle (a_1,a_2),(b_1,b_2);c_1,c_2\,|\,a_1c_2^{-1},a_2c_1^{-1},c_1b_2^{-1}, c_2b_1^{-1}\rangle\\
&&=
\langle (a_1,a_2),(b_1,b_2);\emptys\,|\,a_2b_2^{-1},a_1b_1^{-1}\rangle= \idp_1\diam\idp_1,
\eeqn
where we have used the first two relations to eliminate the internal generators $c_1$ and $c_2$. 

In order to see the naturality of the braiding (axioms \hrel[b2]{b2-3}), i.e. that for any presentation $P:n\to m$
$$
\gammap_{1,m}\circ (\idp_1\otimes P)=(P\otimes \idp_1)\circ \gamma_{1,n}\;\text{and}\;
\gamma_{m,1}\circ(P\otimes \idp_1)  = (\idp_1\otimes P)\circ \gamma_{n,1},
$$
we use the relations appearing in the braiding morphism to eliminate the internal generators in the image/source of $P$.

The  left-  and the right-hand side  of the coassociativity axiom \hrel{a3} are the following:
\beqn
&(\idp_1\diam\Deltap)\circ\Deltap&=\langle (a),(b_1,b_2,b_3);c_1,c_2\,|\,ac_1^{-1},ac_2^{-1},c_1b_1^{-1}, c_2b_2^{-1}, c_2b_3^{-1}\rangle\\
&&=\langle (a),(b_1,b_2,b_3);\emptys\,|\,ab_1^{-1},ab_2^{-1},ab_3^{-1}\rangle\\
&&=\langle (a),(b_1,b_2,b_3);c_1,c_2\,|\,ac_1^{-1},ac_2^{-1},c_1b_1^{-1}, c_1b_2^{-1}, c_2b_3^{-1}\rangle=(\Deltap\diam\idp_1)\circ\Deltap,
\eeqn
where in the first step we use the first two relations to eliminate the internal generators $c_1$ and $c_2$, substituting them elsewhere by $a$, while in the last step we apply the inverse of such a move.

 Analogously, we prove the associativity of the multiplication (axiom \hrel{a1}) by eliminating the internal generators $c_1$ and $c_2$ as follows:
\beqn
&\multp\circ (\idp_1\diam \multp)&=\langle (a_1,a_2,a_3), (b);c_1,c_2\,|\,a_1c_1^{-1},a_2a_3c_2^{-1},c_1c_2b^{-1}\rangle\\
&&=\langle (a_1,a_2,a_3), (b);\emptys\,|\,a_1a_2a_3b^{-1}\rangle \\
&&=\langle (a_1,a_2,a_3), (b);c_1,c_2\,|\,a_1a_2c_1^{-1},a_3c_2^{-1},c_1c_2b^{-1}\rangle=\multp\circ (\multp\diam \idp_1).
\eeqn
 The fact  that $\epsp$ and $\etap$ are respectively left counit and unit (relations \hrel{a4} and \hrel{a2}) is proved below, where we use  the first two relators to eliminate the internal generators $c_1$ and $c_2$. The proofs of \hrel{a4'} and \hrel{a2'} are analogous.
$$
(\epsp\diam\idp_1)\circ \Deltap=\langle (a),(b);c_1,c_2\,|\,ac_1^{-1},ac_2^{-1},c_2b^{-1}\rangle=\langle (a),(b);\emptys\,|\,ab^{-1}\rangle=\idp_1
$$
%\vskip-20pt
$$
\multp\circ(\etap\diam\idp_1)=\langle (a),(b);c_1,c_2\,|\,c_1,ac_2^{-1},c_1c_2b^{-1}\rangle=\langle (a),(b);\emptys\,|\,ab^{-1}\rangle=\idp_1
$$

The bialgebra axiom \hrel{a5} is obtained as follows: \beqn
&&(\multp\diam \multp)\circ (\idp_1\diam\gammap_{1,1}\diam\idp_1)\circ(\Deltap\diam\Deltap)=\\
&&\qquad\quad =
\langle (a_1,a_2), (b_1,b_2);c_1,c_2,c_3,c_4\,|\,a_1c_1^{-1},a_1c_2^{-1},a_2c_3^{-1},a_2c_4^{-1},c_1c_3b_1^{-1},c_2c_4b_2^{-1}\rangle\\
&&\qquad\quad =
\langle (a_1,a_2), (b_1,b_2);\emptys\,|\,a_1a_2b_1^{-1},a_1a_2b_2^{-1}\rangle\\
&&\qquad\quad =
\langle (a_1,a_2), (b_1,b_2);c\,|\,a_1a_2c^{-1},cb_1^{-1},cb_2^{-1}\rangle =
\Deltap\circ\multp,
\eeqn
where in the second step the first four relators are used to eliminate the internal generators, and in the third step  the generator $c$ and the relator $a_1a_2c^{-1}$ is inserted by applying the inverse of move \hrel{ac7}.

Axioms \hrel{a6} and \hrel{a7} follow directly from \hrel{ac7}:
\begin{gather*}
\epsp\circ\multp=\langle (a_1,a_2), \emptys;b\,|\,a_1a_2b^{-1}\rangle=\langle (a_1,a_2), \emptys;\emptys\,|\,\emptys\rangle=\epsp\diam\epsp\\
\Deltap\circ\etap=\langle  \emptys,(a_1,a_2);b\,|\,b, ba_1^{-1},ba_2^{-1}\rangle=\langle \emptys,(a_1,a_2);\emptys\,|\,a_1^{-1},a_2^{-1}\rangle=\etap\diam\etap,
\end{gather*}
\, while axiom \hrel{a8}, i.e. that $\epsp\circ\etap=\idp_0$,  follows from \hrel{ac1}. Finally, the cocommutativity axiom \hrel{a9} holds because the set of relators in a presentation is  unordered.

In order to complete the proof that $1$ is a unimodular bialgebra in the symmetric category $\pres$, it remains to check the integral axioms. \hrel{i2-2'} and \hrel{i3} are obtained by direct application of \hrel{ac1}. The proof of \hrel{i1} is presented below:
\beqn
&(\idp_1\diam\lp)\circ\Deltap&=\langle (a),(b);c_1,c_2\,|\,ac_1^{-1},ac_2^{-1},c_1b^{-1},c_2\rangle=\langle (a),(b);\emptys\,|\,ab^{-1},a\rangle\\
&&=\langle (a),(b);\emptys\,|\,ba^{-1},a\rangle=\langle (a),(b);\emptys\,|\,b,a\rangle=\etap\circ\lp,
\eeqn
where in the first step we have used \hrel{ac7} to eliminate $c_1$ and $c_2$, and in the second and the third steps we have used \hrel{ac5} and \hrel{ac3'}. The proof of \hrel{i1'} is completely analogous.
Then Proposition \ref{antipode/thm} implies that  $1\in {\rm Obj}(\pres)$ is a Hopf algebra with antipode: 
\beqn
&S&=(\idp_1\diam (\lambda\circ\m))\circ(\gammap\diam\idp_1)\circ(\idp_1 \diam(\Deltap\circ \Lp))=\langle (a),(b);c_1,c_2\,|\,c_1b^{-1},c_1c_2^{-1},c_2a\rangle\\
&&=\langle (a),(b);c_2\,|\,bc_2^{-1},c_2a\rangle =\langle (a),(b);\emptys\,|\,ba\rangle=\langle (a),(b);\emptys\,|\,ab\rangle,
\eeqn
where the final equality uses \hrel{ac4}
\end{proof}

\begin{proposition}\label{omega/thm}
There exists a braided monoidal functor $\Psi:\calg\to \pres$ which sends $H$ to 1 and each of the elementary morphisms of $\calg$ to the corresponding morphism of $\pres$ carrying the same name. Moreover for any $n\geq 1$ and $\sigma\in\perm{n}$ 
\begin{itemize}
\item[(a)] $\Psi(\Delta^n)=\langle (a),(b_1,b_2,\dots, b_{n+1});\emptys\,|\,ab_1^{-1},ab_2^{-1},\dots, ab_{n+1}^{-1}\rangle$;
\item[(b)] $\Psi(\mu^n)=\langle (a_1,a_2,\dots, a_{n+1}), (b);\emptys\,|\,a_1a_2\dots a_{n+1}b^{-1}\rangle$;
\item[(c)] $\Psi(\Upsilon^\sigma)=\langle (a_1,a_2,\dots, a_n),(b_1,b_2,\dots, b_n);\emptys\,|\,a_1b_{\sigma(1)}^{-1},a_2b_{\sigma(2)}^{-1},\dots, a_nb_{\sigma(n)}^{-1}\rangle$;
\item[(d)] $\Psi(\lambda\circ\eta)=\langle \emptys, \emptys;\emptys\,|\, 1\rangle$.
\end{itemize}
\end{proposition}

\begin{proof}
The existence of a braided monoidal functor $\Psi:\calg\to \pres$ follows from Theorem \ref{Hpres/thm} and the universal property of $\calg$. 

{\sl (a)} follows by induction on $n$. Indeed, for $n=1$ it is true by the definition of $\Psi$ and if the statement holds for $\Psi(\Delta^{n-1})$, then
\begin{eqnarray*}
&\Psi(\Delta^{n})& =\Psi(\Delta^{n-1}\diam\id_1)\circ\Psi(\Delta)\\&&=
\langle (a),(b_1,\dots, b_{n}, b_{n+1});c_1,c_2\,|\,c_1b_1^{-1},c_1b_2^{-1},\dots, c_1b_{n}^{-1},c_2b_{n+1}^{-1},ac_1^{-1},ac_2^{-1}\rangle.
\end{eqnarray*}
Now it is enough to use the last two relators to eliminate $c_1$ and $c_2$, substituting them elsewhere by $a$. The proof of (b) is similar and it is left to the reader.

{\sl(c)} follows by induction on the minimal number $m_\sigma$ of transpositions of the type $(i, i+1)$, $1\leq i<n$ in which $\sigma$ decomposes. Indeed, the statement for $\sigma=(i, i+1)$, reduces to the definition of $\Psi(\braid)$. On the other hand,  Theorem \ref{symm-coherence/thm} (iii) implies that, if the statement holds for $\sigma$ and $\sigma'$, then:
\begin{eqnarray*}
&\Psi (\Upsilon^{\sigma'\circ\sigma})&=\Psi (\Upsilon^{\sigma'})\circ \Psi (\Upsilon^{\sigma})\\
&&= \langle (a_1,\dots, a_n),(b_1,\dots, b_n);c_1,\dots,c_n\,|\,a_1c_{\sigma(1)}^{-1},\dots, a_nc_{\sigma(n)}^{-1},c_1b_{\sigma'(1)}^{-1},\dots, c_nb_{\sigma'(n)}^{-1}\rangle\\
&&= \langle (a_1,\dots, a_n),(b_1,\dots, b_n); \emptys\,|\,a_{\sigma^{-1}(1)}b_{\sigma'(1)}^{-1},\dots, a_{\sigma^{-1}(n)}b_{\sigma'(n)}^{-1}\rangle\\
&&= \langle (a_1,\dots, a_n),(b_1,\dots, b_n); \emptys\,|\,a_{1}b_{\sigma'\circ\sigma(1)}^{-1},\dots, a_{n}b_{\sigma'\circ\sigma(n)}^{-1}\rangle.
\end{eqnarray*}

{\sl(d)} is obtained as follows:
$$\Psi(\lambda\circ\eta)=\langle \emptys, \emptys;a\,|\, a,a\rangle
=\langle \emptys, \emptys;a\,|\, 1,a\rangle
=\langle \emptys, \emptys;\emptys\,|\, 1\rangle$$
\end{proof}
 
\begin{figure}[htb]
  \includegraphics{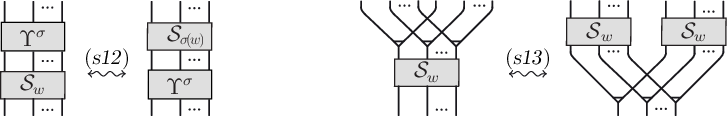}
\caption{Properties of $S_w$.}
\label{Sw-prop/fig}\label{E:s12} \label{E:s13}
%\vskip-6pt
\end{figure}

\begin{lemma}\label{Sw/thm} To any  word $w$ of length $l\geq 0$ in $G=\{g_1,g_2,\dots,g_n\}$ we associate the morphism $\cS_w:H^{l}\to H^{l}$ in $\calg$, defined inductively on $l$ as follows:
$$
\cS_{1}=\id_\one,\qquad \cS_{g_i}=\id_1,\qquad \cS_{g_i^{-1}}=S,\qquad \cS_{w_1w_2}=\cS_{w_1}\diam\cS_{w_2}.
$$
 Then (see Proposition \ref{delta-mu-n/thm} for the definition of $\pi_l$)
\begin{gather*}\Upsilon^\sigma\circ S_{w}=S_{\sigma(w)}\circ \Upsilon^\sigma, \sigma\in\perm{l}, \tag*{\hrel{s12}}
\\
\Upsilon^{\pi_l}\circ\Delta^{\diam l}\circ S_w=(S_w\diam S_w)\circ \Upsilon^{\pi_l}\circ\Delta^{\diam l}.  \tag*{\hrel{s13}}
\end{gather*}
\end{lemma}
The statement of  Lemma \ref{Sw/thm} is illustrated in Figure \ref{Sw-prop/fig}.

\begin{proof}
\hrel{s12} follows from the naturality of $\Upsilon$, while \hrel{s13} is obtained from \hrel{s5} in Table \ref{antipode-form/fig} using induction on $l$. 
\end{proof}

We can now complete the proof of Theorem \ref{main1/thm}.

\begin{theorem}\label{bomega/thm}
There exists a braided monoidal functor $\overline{\Psi}:\pres\to \calg$ such that $\overline{\Psi}\circ\Psi=\id_{\calg}$ and $\Psi\circ\overline{\Psi}=\id_{\pres}$. In particular, $\Psi$ and $\overline{\Psi}$ are category equivalences.
\end{theorem}

\begin{proof}
On the set of objects we define
$\overline{\Psi}(n)=H^{n}$.
 Then, given a morphism in $\pres$
$$
P=\langle\, (a_1,\dots,a_n),\,(b_1,\dots,b_m);\, C\, |\, R\, \rangle:n\to m,
$$
  the definition of $\overline{\Psi}(P)$ involves some choices. 
 
 First of all we fix a total order on the set of internal generators $C=(c_1,\dots,c_k)$ and denote by $\bar G$ the  ordered set of generators of $P$, starting with the ones in the source, then the  internal ones and finally the ones in the target:
$$
\bar G=(g_1,g_2,\dots, g_{r})\defeq(a_1,\dots,a_n,c_1,\dots,c_k,b_1,\dots,b_m),\; r=n+m+k. $$
Then we fix a total order on the set of relators and denote this ordered set by  $\bar R=(w_1,\dots,w_s)$.  Let  
$$
l_i=|w_i|, \quad w_P=w_1w_2\dots w_s,\quad e_i=|w_P|_{g_i}.
$$
 Let also $w^+_P$ be the word obtained from $w_P$ by inverting the signs of all negative exponents of the generators and
 $$
 l=|w_P^+|=|w_P|=\sum_{i=1}^s l_i=\sum_{i=1}^r e_i.$$
Finally, we choose $\sigma_P\in \perm{l}$ such that $\sigma_P(\prod_{i=1}^r g_i^{e_i})=w^+_P$ (see p.\pageref{action-sym/par}) and define:
$$
\overline{\Psi}(P)=\left(\left(\diam_{i=1}^s (\lambda\circ \mu^{l_i-1}\circ\cS_{w_i})\circ
\Upsilon^{\sigma_P}\circ
(\diam_{j=1}^r \Delta^{e_j-1})\right)\diam\id_{m}\right)\circ
(\id_{n}\diam \Lambda^{\diam k}\diam\coev_m),
$$
where $\Lambda^{\otimes 0}:=\id_\one$. The definition of $\overline{\Psi}$ is illustrated in Figure \ref{baromega/fig} and a specific example is shown in Figure \ref{example1/fig}.

\begin{figure}[htb]
  \includegraphics{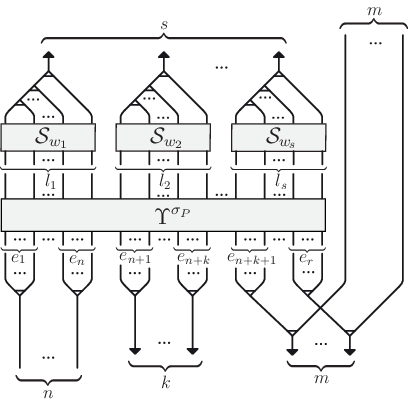}
\caption{Definition of $\overline{\Psi}(P)$}
\label{baromega/fig}
%\vskip-6pt
\end{figure} 

\begin{figure}[htb]
  \includegraphics{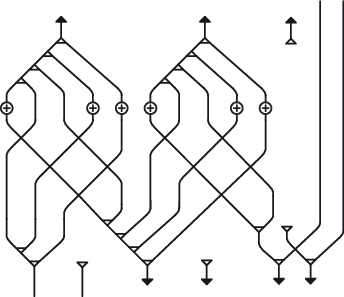}
\caption{$\overline{\Psi}(\langle (a_1,a_2),(b_1,b_2);c_1,c_2\,|\, c_1^{-1}a_1c_1a_1^{-2},b_1^{-1}c_1b_1c_1^{-2},1\rangle)$}
\label{example1/fig}
%\vskip-6pt
\end{figure} 

\medskip
\begin{lemma}\label{indipchoices/thm}
$\overline{\Psi}(P)$ is independent of the choice of $\sigma_P$ and of the choices of total order of the sets of internal generators and relators. 
\end{lemma}

\begin{proof}
Let $\sigma_P$ and $\sigma'_P$ be two permutations in $\perm{l}$ that rearrange the word $\prod_{i=1}^r g_i^{e_i}$ into the word $w^+_P$. Then  $\sigma_P^{-1}\circ\sigma_P'$ leaves unchanged   $\prod_{i=1}^r g_i^{e_i}$; hence $\sigma_P'=\sigma_P\circ (\sigma_1\times\sigma_2\times\dots\times\sigma_r)$ for some permutations $\sigma_i\in \perm{e_i}$, $i=1,\dots,r$ (Definition \ref{iota/defn}) and according to Theorem \ref{symm-coherence/thm} (iii) and (iv) we have:
$$
\Upsilon^{\sigma'_P}=\Upsilon^{\sigma_P}\circ (\Upsilon^{\sigma_1}\diam\Upsilon^{\sigma_2}\diam\dots\diam \Upsilon^{\sigma_r}).
$$
 Then \hrel{a11} in Proposition \ref{delta-mu-n/thm} implies that (see Figure \ref{indepsigma/fig})
$$
\Upsilon^{\sigma'_P}\circ
\left(\diam_{j=1}^r \Delta^{e_j-1}\right)=
\Upsilon^{\sigma_P}\circ\left(\diam_{j=1}^r (\Upsilon^{\sigma_j}\circ\Delta^{e_j-1}) \right)=
\Upsilon^{\sigma_P}\circ\left(\diam_{j=1}^r \Delta^{e_j-1}\right).
$$
Therefore substituting $\sigma_P$ with $\sigma_P'$ doesn't change $\overline{\Psi}(P)$. 

\begin{figure}[htb]
  \includegraphics{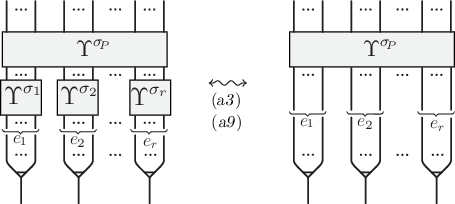}
\caption{Independence of the choice of $\sigma_P$.}
\label{indepsigma/fig}
%\vskip-6pt
\end{figure} 

In order to prove the independence of the total order on the set of relators, it is enough to analyze what happens if  two consecutive relators $r_i$ and $r_{i+1}$ are interchanged. Such a move changes the permutation $\sigma_P$ by a composition with a permutation which interchanges the positions of the $l_i$-tuple and the $l_{i+1}$-tuple (Figure \ref{indepR/fig}). Therefore in the expression for $\overline{\Psi}(P)$, $\Upsilon^{\sigma_p}$ is substituted by $(\id_{l'_i}\diam \braid_{l_i,l_{i+1}}\diam \id_{l''_i})\circ \Upsilon^{\sigma_p}$, where $l'_i=\sum_{j=1}^{i-1}l_j$ and $l''_i=\sum_{j=i+2}^{s}l_j$. Then by the naturality of $\braid$ (axiom \hrel{b3} in Table \ref{axiomsH/fig}) we have:
$$
\left((\lambda\circ \mu^{l_{i+1}-1}\circ\cS_{w_{i+1}})\diam (\lambda\circ \mu^{l_{i}-1}\circ\cS_{w_{i}})\right)\circ
\braid_{l_i,l_{i+1}}= \braid_{\one ,\one }\circ \left((\lambda\circ \mu^{l_{i}-1}\circ\cS_{w_{i}})\diam (\lambda\circ \mu^{l_{i+1}-1}\circ\cS_{w_{i+1}})\right).
$$
Hence interchanging $r_i$ and $r_{i+1}$ doesn't change $\overline{\Psi}(P)$.
\begin{figure}[htb]
  \includegraphics{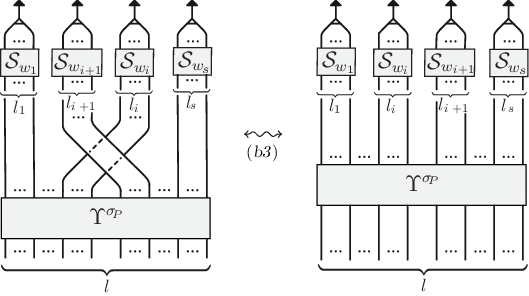}
\caption{Independence of the choice of total order on the set of relators.}
\label{indepR/fig}
%\vskip-6pt
\end{figure} 

The proof of the independence of the total order on the set of internal generators is completely analogous and it is illustrated in Figure \ref{indepG/fig}.
\begin{figure}[htb]
  \includegraphics{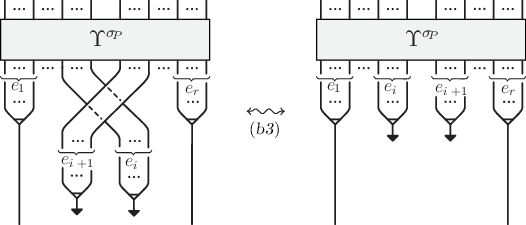}
\caption{Independence of the choice of total order on the set of internal generators.}
\label{indepG/fig}
%\vskip-6pt
\end{figure} 
\end{proof}

\begin{lemma}\label{ACinv/thm}
$\overline{\Psi}(P)$ depends only on the AC-equivalence class of  $P$.
\end{lemma}

\begin{proof} We have to show that  $\overline{\Psi}(P)$ is invariant under the five AC-moves in Definition \ref{relativepr/def}.

\medskip
\noindent{\sl Invariance under \hrel{ac1}.} Let $Q$ be the presentation obtained from $P$ by adding a new internal generator $c_0$ and a new relator $w_0=c_0v_0$, where $v_0$ is any word in $G$:
$$
Q=\langle\, (a_1,\dots,a_n),\,(b_1,\dots,b_m);\, c_0,c_1,\dots,c_k\, |\, c_0v_0,w_1,\dots,w_s\, \rangle
$$
Let  
$$l_0=|w_0|,\;l'=l+l_0,\;\; f_i=|w_0|_{g_i},\;  e_0'=1 \;  \text{and}\; e_i'=|w_Q^+|_{g_i}=e_i+f_i \;\text{for}\; 1\leq i\leq r.$$ 
  Assuming that the new generator and relator are the first ones, we have 
$$
\overline{\Psi}(Q)=\left(\left(\diam_{i=0}^s (\lambda\circ \mu^{l_i-1}\circ\cS_{w_i})\circ
\Upsilon^{\sigma_{Q}}\circ
\Delta^\bfe\right)\diam\id_{m}\right)\circ
(\id_{n}\diam \Lambda^{\diam (k+1)} \diam\coev_m),
$$
where 
$$\Delta^\bfe=(\diam_{j=1}^n \Delta^{e'_j-1})\otimes \id_1\otimes (\diam_{j=n+1}^r \Delta^{e'_j-1}).$$
Moreover, since $\sigma_Q\in \perm{l'}$ is defined up to multiplication on the right by element in $\times_{i=1}^r\perm{e_i'}$, we can choose it in such a way that (see Figure \ref{inv-ac1/fig}):
$$
\Upsilon^{\sigma_Q}=(\id_1\otimes \Upsilon^{\sigma})\circ (\braid_{f, 1}\otimes\id_{l'-1-f}),
$$
where $f=\sum_{i=1}^ne_i'$ and $\sigma\in\perm{l'-1}$ is such that $\sigma(\prod_{i=1}^r g_i^{e_i'})=v_0^+\,w_P^+$. 
Then, using the naturality of $\Upsilon$,  $\overline{\Psi}(Q)$ can be written as (Figure \ref{inv-ac1/fig}):
\beqn
&\overline{\Psi}(Q)&=
\left(\left((\lambda\circ \mu^{l_0-1})\circ (\Lambda\diam \cS_{v_0})\right)\diam 
\left(\diam_{i=1}^s (\lambda\circ \mu^{l_i-1}\circ\cS_{w_i})\right)\diam\id_{m}\right)\\
&&\quad\;\circ
\left( (\Upsilon^{\sigma}\circ(\diam_{j=1}^r \Delta^{e_j'-1}))\diam\id_{m}\right)
\circ
(\id_{n}\diam \Lambda^{\diam k}\diam\coev_m)
\eeqn
Now applying the integral properties \hrel{i2}, \hrel{i3} and the invariance of the counit under the action of the antipode \hrel{s7}, we have (see Figure \ref{inv-ac12/fig}):
\beqn
&\overline{\Psi}(Q)&=
\left(\diam_{i=1}^s (\lambda\circ \mu^{l_i-1}\circ\cS_{w_i})\diam\id_m\right)\\
&&\quad\;\circ
\left( \left( \left(\epsilon^{\diam (l_0-1)}\diam\id_l\right)\circ\Upsilon^{\sigma}\right)\circ(\diam_{j=1}^r \Delta^{e_j'-1}))\diam\id_{m}\right)
\circ
(\id_{n}\diam \Lambda^{\diam k}\diam\coev_m)
\eeqn
Finally, using again the naturality of $\Upsilon$ and the counit property \hrel{a4}  we obtain:
$$\left(\left(\epsilon^{\diam (l_0-1)}\diam\id_l\right)\circ\Upsilon^{\sigma}\right)\circ(\diam_{j=1}^r \Delta^{e_j'-1})=\Upsilon^{\sigma_P}\circ(\diam_{j=1}^r ((\epsilon^{\diam f_j}\diam\id_{e_j})\circ  \Delta^{e_j'-1}))=\Upsilon^{\sigma_P}\circ (\diam_{j=1}^r \Delta^{e_j-1}),
$$
 where as desired $\sigma_P(\prod_{i=1}^r g_i^{e_i})=w_P^+$.
 
\begin{figure}[htb]
  \includegraphics{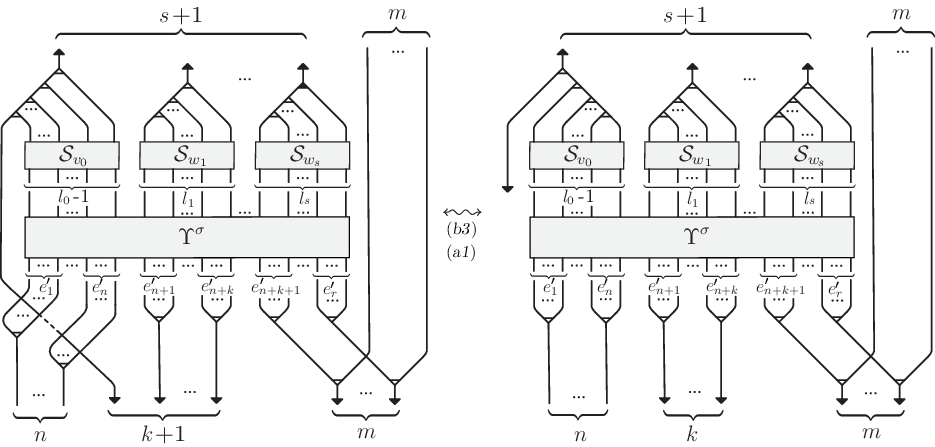}
\caption{Invariance of $\overline{\Psi}(P)$ under \hrel{ac1}, part I.}
\label{inv-ac1/fig}
%\vskip-6pt
\end{figure}

\begin{figure}[hbt]
  \includegraphics{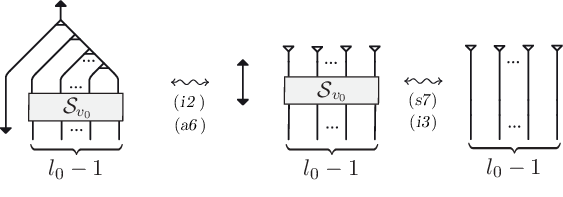}
\caption{Invariance of $\overline{\Psi}(P)$ under \hrel{ac1}, part II.}
\label{inv-ac12/fig}
%\vskip-6pt
\end{figure} 

\smallskip
\noindent{\sl Invariance under \hrel{ac2}.} The move \hrel{ac2} consists of  multiplying  a relator in $R$ on the left by $g_i g_i^{-1}$ or by $g_i^{-1}g_i$ for some  $g_i\in G$. Without loss of generality we can assume  that the relator  is the first one. The resulting presentation is:
$$
Q=\langle\, (g_1,\dots,g_n),\,(g_{n+k+1},\dots,g_r);\, g_{n+1},\dots,g_{n+k}\, |\, w_1',\dots,w_s'\, \rangle,
$$
where $w_1'=g_i g_i^{-1}w_1$ and $w_j'=w_j$, $1<j\leq s$. 
Let  
$$e_h'=|w_Q|_{g_h},\,1\leq h\leq r,\quad l_j'=|w_j'|,\,1\leq j\leq s,\quad f= \sum_{j=1}^{i-1}e'_j.
$$
Then
$$
\overline{\Psi}(Q)=\left(\left(\diam_{i=1}^s (\lambda\circ \mu^{l_i'-1}\circ\cS_{w_i'})\circ
\Upsilon^{\sigma_{Q}}\circ
(\diam_{j=1}^r \Delta^{e_j'-1})\right)\diam\id_{m}\right)\circ
(\id_{n}\diam \Lambda^{\diam k} \diam\coev_m),
$$
where we choose $\sigma_Q\in\perm{l+2}$ in the form
$\sigma_Q=(\id_{\perm{2}}\times\sigma_P)\circ\sigma$ with
$$
\sigma(j)=\begin{cases} j+2\;\;\text{if}\; 1\leq j\leq f,\\
1\;\text{if}\;\; j= f+1,\\
2\;\text{if}\; \;j= f+2,\\
j\;\text{if}\; \; f+3\leq j\leq l+2.\end{cases}
$$
The relevant part of $\overline{\Psi}(Q)$ is presented on the left-hand side of Figure \ref{inv-ac2/fig}. As indicated in the figure, now the equivalence of $\overline{\Psi}(Q)$ and $\overline{\Psi}(P)$ follows from the naturality of $\braid$  and the defining axioms of the antipode, the unit and the counit. The proof of the invariance under  multiplication by $g_i^{-1}g_i$ is analogous.

\begin{figure}[htb]
  \includegraphics{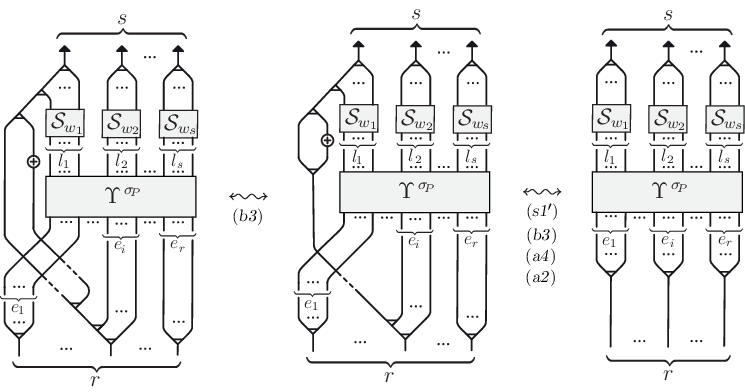}
\caption{Invariance of $\overline{\Psi}(P)$ under \hrel{ac2}.}
\label{inv-ac2/fig}
%\vskip-6pt
\end{figure} 

\smallskip
\noindent{\sl Invariance under \hrel{ac3}.}  
Since $\overline{\Psi}$ doesn't depend on the order of the relators, it is enough to show that  $\overline{\Psi}(P)=\overline{\Psi}(Q)$, where $Q$ is obtained by multiplying the second relator by the first:
$$
Q=\langle\, (g_1,\dots,g_n),\,(g_{n+k+1},\dots,g_r);\, g_{n+1},\dots,g_{n+k}\,|\; w_1,w_1w_2,\dots,w_s\, \rangle
$$
If $w_1=1$, the move does not change the presentation, hence assume $l_1\geq1$. Then (see Figure \ref{inv-ac31/fig}):
\beqn
&&(\lambda\circ\mu^{l_1-1}\circ\S_{w_1})\diam(\lambda\circ\mu^{l_2-1}\circ\S_{w_2})\\
&&\quad=\lambda\circ\mu\circ (((\lambda\diam\id_1)\circ(\Delta\circ\mu^{l_1-1}\circ\S_{w_1}))\diam (\mu^{l_2-1}\circ\S_{w_2}))\\
&&\quad= (\lambda\diam(\lambda\circ\mu))\circ(((\mu^{l_1-1}\diam\mu^{l_1-1})\circ\Upsilon^{\pi_{l_1}}\circ \Delta^{\diam l_1}\circ\S_{w_1})\diam(\mu^{l_2-1}\circ\S_{w_2}))\\
&&\quad=(\lambda\diam(\lambda\circ\mu))\circ((\mu^{l_1-1}\circ\S_{w_1}) \diam(\mu^{l_1-1}\circ\S_{w_1})\diam(\mu^{l_2-1}\circ\S_{w_2}))\circ((\Upsilon^{\pi_{l_1}}\circ \Delta^{\diam l_1})\diam\id_{l_2})\\
&&\quad=((\lambda\circ\mu^{l_1-1}\circ\S_{w_1}) \diam(\lambda\circ\mu^{l_1+l_2-1}\circ\S_{w_1w_2}))\circ((\Upsilon^{\pi_{l_1}}\circ \Delta^{\diam l_1})\diam\id_{l_2})
\eeqn
After substituting the above expression in $\overline{\Psi}(P)$ and applying the naturality of  $\Upsilon$  we obtain (see Figure \ref{inv-ac3/fig}):
$$
((\Upsilon^{\pi_{l_1}}\circ \Delta^{\diam l_1})\diam\id_{l-l_1})\circ\Upsilon^{\sigma_P}\circ (\diam_{i=1}^r\Delta^{e_i-1})=
\Upsilon^{\sigma_Q}\circ (\diam_{i=1}^r\Delta^{e_i'-1}),
$$
where $e_i'=e_i+|w_1|_{g_i}$. Therefore $\overline{\Psi}(P)=\overline{\Psi}(Q)$.

\begin{figure}[htb]
  \includegraphics{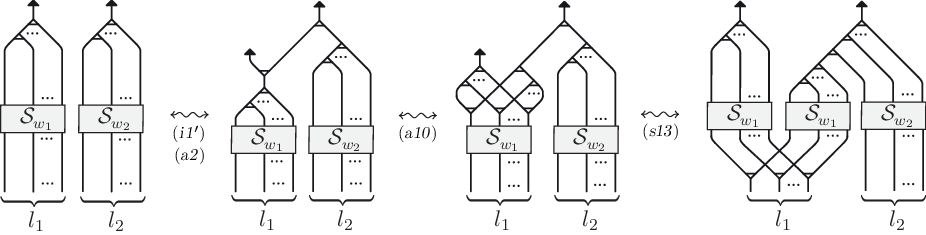}
\caption{Invariance of $\overline{\Psi}(P)$ under \hrel{ac3} -- part I.}
\label{inv-ac31/fig}
%\vskip-6pt
\end{figure} 

\begin{figure}[htb]
  \includegraphics{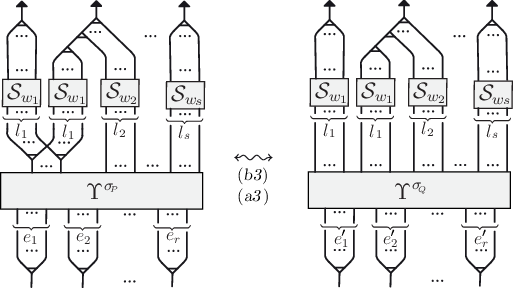}
\caption{Invariance of $\overline{\Psi}(P)$ under \hrel{ac3} -- part II.}
\label{inv-ac3/fig}
%\vskip-6pt
\end{figure} 

\noindent{\sl Invariance under \hrel{ac4}.} 
 The move \hrel{ac4} permutes cyclically  the letters of a relator $w=g_{i_1}g_{i_2}\dots g_{i_p}$ transforming it into  $w'=g_{i_p}g_{i_1}g_{i_2}\dots g_{i_{p-1}}$, $p\geq 1$. In the expression for $\overline{\Psi}(P)$ this corresponds to the replacement of $\lambda\circ \mu^{p-1}\circ\cS_{w}$ by $\lambda\circ \mu^{p-1}\circ\cS_{w'}\circ \braid_{p-1,1}$, but as it is proved below the last two morphisms are equal in $\calg$ (cf. Figure \ref{inv-ac4/fig}): 
\beqn
&\lambda\circ\mu^{p-1}\circ\cS_{w'}\circ \braid_{p-1,1}&=\lambda\circ\mu\circ(\id_1\diam\mu^{p-2}) \circ \braid_{p-1,1}\circ\cS_{w}=\lambda\circ\mu\circ\braid_{1,1}\circ (\mu^{p-2}\diam\id_1)\circ\cS_{w}\\
&&=\lambda\circ\mu^{p-1}\circ\cS_{w},
\eeqn
where  we first apply \hrel{s12} in Lemma \ref{Sw/thm}, then the naturality of $\braid$ and the associativity of $\mu$, and finally we apply \hrel{i6} in Proposition \ref{antipode/thm}. 

\begin{figure}[htb]
  \includegraphics{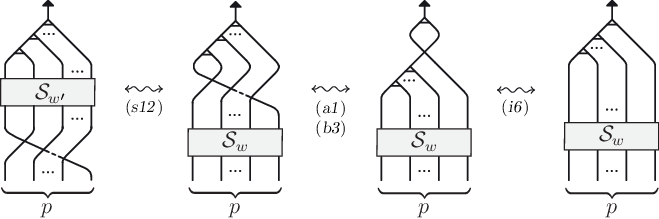}
\caption{Invariance of $\overline{\Psi}(P)$ under \hrel{ac4}.}
\label{inv-ac4/fig}
%\vskip-6pt
\end{figure} 

\smallskip
\noindent{\sl Invariance under \hrel{ac5}.} 
Let $w=g_{i_1}^{\epsilon_1}g_{i_2}^{\epsilon_2}\dots g_{i_p}^{\epsilon_p}$ be the relator that is replaced by its inverse. In the expression for $\overline{\Psi}(P)$ replacing  $w$ by $w^{-1}$ corresponds to the replacement of $\lambda\circ \mu^{p-1}\circ\cS_{w}$ by $\lambda\circ \mu^{p-1}\circ\cS_{w^{-1}}\circ \Upsilon^{\tau_p}$ (cf. Proposition \ref{delta-mu-n/thm}). As it is proved below the last two morphisms are equal in $\calg$ (cf. Figure \ref{inv-ac5/fig}): 
$$\lambda\circ\mu^{p-1}\circ\cS_{w}=\lambda\circ S\circ\mu^{p-1}\circ\cS_{w}
=\lambda\circ\mu^{p-1}\circ \Upsilon^{\tau_p}\circ S^{\diam p}\circ\cS_{w}
=\lambda\circ\mu^{p-1}\circ \cS_{w^{-1}}\circ \Upsilon^{\tau_p},
$$
where  we first apply \hrel{i5}, then \hrel{s11} in Proposition \ref{delta-mu-n/thm}, and finally \hrel{s12} in Lemma \ref{Sw/thm} and \hrel{s2} in Proposition \ref{antipode/thm}. 
\end{proof}

\begin{figure}[htb]
  \includegraphics{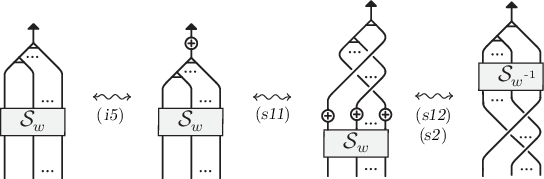}
\caption{Invariance of $\overline{\Psi}(P)$ under \hrel{ac5}.}
\label{inv-ac5/fig}
%\vskip-6pt
\end{figure} 

We can now complete the proof of Theorem \ref{bomega/thm}. We will first show that $\overline{\Psi}$ defines a braided monoidal functor from $\pres$ to $\calg$, i.e. that it preserves identities, compositions, products and braiding structure. 

Using the definition of $\idp_n$ and those of the evaluation and coevaluation morphisms in
Proposition \ref{lambda/thm}, we obtain (see Figure \ref{bomegaid/fig}):
\beqn
&\overline{\Psi}(\idp_n)&=\left(\left( (\lambda\circ \mu\circ (\id_1\diam S))^{\diam n}\circ
\Upsilon^{\pi_n^{-1}}\right)\diam\id_{n}\right)\circ (\id_{n}\diam \coev_n)\\
&&=(\ev_n\diam\id_n )\circ (\id_n\diam\coev_n)=\id_n.
\eeqn

\begin{figure}[htb]
  \includegraphics{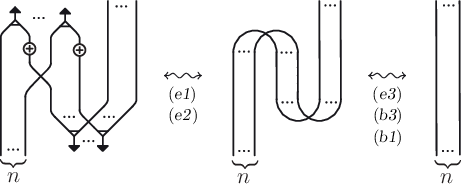}
\caption{$\overline{\Psi}(\idp_n)=\id_n$.}
\label{bomegaid/fig}
%\vskip-6pt
\end{figure} 

To see that $\overline{\Psi}$ preserves compositions, let 
\begin{gather*}
P=\langle\, \overline{A}, (h_1,\dots,h_m);\, C\, |\, R\, \rangle:n\to m,\quad 
Q=\langle\, (h_1,\dots,h_m),\,\overline{B};\, C'\, |\, R'\, \rangle:m\to n',  
\\
Q\circ P=\langle\, \overline{A}, \overline{B};\{h_1,\dots,h_m\}\sqcup C\sqcup C'\, |\, R\sqcup R'\, \rangle:n\to n'\, ,
\\
e_i=|w^+_P|_{h_i}, \; e_i'=|w^+_Q|_{h_i}, \;1\leq i\leq m,\quad p=|w^+_P|-\sum_{i=1}^{m} e_{i} ,\quad q=|w^+_Q|-\sum_{i=1}^m e'_{i}.
\end{gather*}
Using the naturality of $\braid$, the coassociativity and the counit axioms (if some of the $e_i$'s or $e_i'$'s are 0), we can transform the essential part of the morphism $\overline{\Psi}(Q)\circ \overline{\Psi}(P)$ as follows (cf. Figure \ref{bomegacomp/fig}):
\beqn
&&(\Upsilon^{\sigma_P}\diam\Upsilon^{\sigma_Q})\circ(\id_p\diam(((\diam_{i=1}^{m} \Delta^{e_{i}-1})\diam(\diam_{i=1}^{m} \Delta^{e'_{i}-1}))\circ\Upsilon^{\pi_m}\circ(\Delta\circ \Lambda)^{\diam m})\diam\id_q)\\
&&\qquad\qquad=\Upsilon^{\sigma_{Q\circ P}}\circ(\id_p\diam (\diam_{i=1}^{m} (\Delta^{e_{i}+e'_{i}-1}\circ \Lambda))\diam\id_q).
\eeqn
Therefore $\overline{\Psi}(Q)\circ \overline{\Psi}(P)=\overline{\Psi}(Q\circ P)$.
 
\begin{figure}[htb]
  \includegraphics{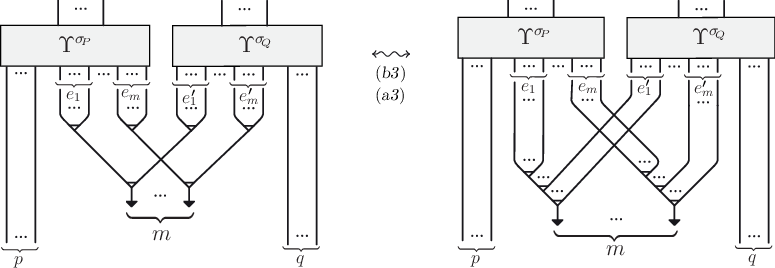}
\caption{$\overline{\Psi}(Q)\circ \overline{\Psi}(P)=\overline{\Psi}(Q\circ P)$.}
\label{bomegacomp/fig}
%\vskip-6pt
\end{figure} 

\smallskip
Since the product in $\pres$ is given by the disjoint union on the set of generators and relators, $\overline{\Psi}(P\diam Q)=\overline{\Psi}(P)\diam \overline{\Psi}(Q)$ follows from the naturality of the braiding, while the proof that $\overline{\Psi}$ preserves the braiding morphism is analogous to the one for the identity and is illustrated in Figure \ref{bomega-gamma/fig}.

\begin{figure}[htb]
  \includegraphics{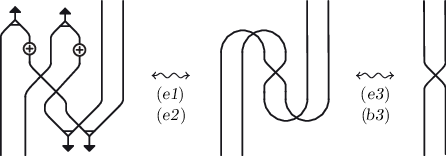}
\caption{$\overline{\Psi}(\dot{\gamma}_{1,1})=\gamma_{1,1}$.}
\label{bomega-gamma/fig}
%\vskip-6pt
\end{figure} 

Now it remains to show that $\overline{\Psi}$ is the inverse functor of $\Psi$, i.e. that $\overline{\Psi}\circ\Psi=\id_{\calg}$ and  $\Psi\circ\overline{\Psi}=\id_\pres$. The first identity  follows from the fact that if $F$ is any of the generating morphisms of $\calg$ ($\braid$, $\mu$, $\Delta$,  $\eta$, $\epsilon$, $\lambda$ and $\Lambda$), we have $\overline{\Psi}(\dot F)=F$. The verification of this is straightforward  and the only nontrivial cases are shown in Figure \ref{composition1/fig}.

\begin{figure}[htb]
  \includegraphics{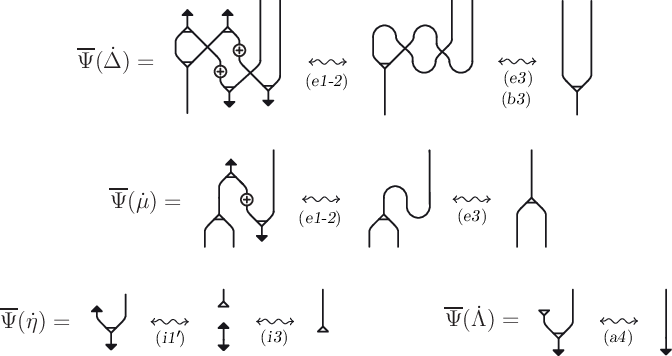}
\caption{Proof that $\overline{\Psi}\circ\Psi=\id_{\calg}$.}
\label{composition1/fig}
%\vskip-6pt
\end{figure} 

To see that $\Psi\circ\overline{\Psi}=\id_\pres$, let
$$
P=\langle (g_1,\dots,g_n), (g_{n+k+1},\dots, g_{n+k+m}); g_{n+1},\dots,g_{n+k}\,|\,w_1,\dots, w_s\rangle: n\to m
$$
be a morphism in $\pres$. The definitions of $\overline{\Psi}$ and $\Psi$  (Proposition \ref{omega/thm}) imply that $\Psi\circ\overline{\Psi}(P)$ is AC-equivalent to:
$$
\langle (g_1,\dots,g_n), (g_{n+k+1},\dots, g_{n+k+m}); g_{n+1},\dots, g_{n+k}, c_1,\dots,c_m\,|\, w_1,\dots, w_s,c_1 g_{n+k+1}^{-1},\dots ,c_m g_{n+k+m}^{-1} \rangle.
$$
Then $\Psi\circ\overline{\Psi}(P)= P$ follows by  using the relators $c_i g_{n+k+i}^{-1}$ to eliminate all $c_i$'s, substituting them by the $g_{n+k+i}$'s.
\end{proof}
%%%%%%%%%%%%%%

%%%%%%%%%%%%%%%%%%%%%
\section{$\calg$ as quotient of the category $\Algf$}
\label{Hr-Htr/sec}
%%%%%%%%%%%%%%%%%%%%%

This section is dedicated to the proof of Theorem \ref{main2/thm}.  We start by reviewing the definition of the category $\RHB$ and refer the reader to \cite{GS99} and \cite[Subsection~4.3 ]{BBDP23} for a detailed discussion.

\begin{figure}[htb]
 \centering
 \includegraphics{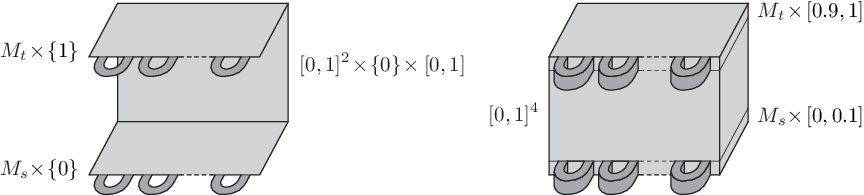}
 \caption{Canonical realization of $M_{s,t}$ and $M_{s,t} \times [0,1]$ in $\R^4$.}
 \label{cobordism01/fig}
\end{figure}

Let $M = M_{s,t} \cong M_s \bcsum M_t$ be the boundary connected sum of two (absolute) connected oriented \dmnsnl{3} \hndlbds{1}
\[
M_s \cong H^0 \cup_{i=1}^s H^1_i 
\text{ \ and \ } 
M_t \cong H^0 \cup_{i=1}^t H^1_i,
\]
with $s,t \geqs 0$. We assume that $M_s$ and $M_t$ are canonically realized inside $\R^3$, by identifying $H^0$ with $[0,1]^3$ and attaching each 1-handle $H^1_i$ to $\left] 0,1 \right[^2 \times \{1\}$. So, we can set 
\[
M_{s,t} = (M_s \times \{0\}) \cup
([0,1]^2 \times \{0\} \times [0,1]) \cup (M_t \times \{1\}) \subset \R^4,
\] 
as depicted on the left-hand side of Figure~\ref{cobordism01/fig}. 

An oriented \dmnsnl{4} {\sl relative \hndlbd{2}} built on $M_{s,t}$ is an oriented smooth \mnfld{4} with a given handle decomposition
\[
W = M_{s,t}\times [0,1] \cup_{i=1}^m H^1_i \cup_{j=1}^n H^2_j,
\]
where $H^1_i$ and $H^2_j$ are 1- and 2-handles attached to the front boundary $M_{s,t}\times\{1\}$ of $M_{s,t}\times [0,1] $. We observe that $M_{s,t} \times [0,1]$ can be identified canonically with a 4-dimensional 1-handlebody with a single 0-handle and \dmnsnl{4} 1-handles as shown on the right in Figure \ref{cobordism01/fig}.

The oriented 4-dimensional relative 2-handlebodies built over the $3$-manifolds $M_{s,t}$ with $s,t \geqs 0$ form a strict monoidal category $\RHB$ whose objects are connected oriented \dmnsnl{3} \hndlbds{1} $M_s$ for $s \geqs 0$, and whose morphisms from $M_s$ to $M_t$ are $2$-equivalence classes of oriented \dmnsnl{4} relative \hndlbds{2} built on $M_{s,t}$, where two morphisms $W,W':M_s\to M_t$ are 2-equivalent if their relative \hndlbd{2} decompositions are related by a \textsl{$2$-deformation}, meaning a finite sequence of the following operations:
\begin{itemize}
\item[\(a)]
isotoping the attaching maps of the handles;
\item[\(b)]
adding/deleting a canceling pair consisting of a 1-handle and a 2-handle;
\item[\(c)]
sliding a 2-handle over another one.
\end{itemize}

%The operation of sliding a 1-handle over another one is also admitted, as it can be obtained from \(b) and \(c), see for instance \cite[Figure~2.2.11]{BP11}. Since the 2-deformations preserve the diffeomorphism type of the handlebody, if two oriented \dmnsnl{4} relative \hndlbds{2} are \qvlnt{2}, then they are diffeomorphic. 

The composition of two morphisms in $\RHB$  is obtained by  taking their vertical juxtaposition, with $W'$ on top of $W$, and gluing the two morphisms, by identifying canonically the target of the first with the source of the second. On the other hand, the monoidal product is given by horizontal juxtaposition and boundary connected sum of both objects and morphisms.

For each $s \geqs 0$, the identity $\id_{M_s}$ is represented by the product $M_s \times [0,1]$ with the natural handlebody decomposition. In particular, $\idone = M_0 \times [0,1]$, since $\one = M_0$.

\smallskip 
Following \cite{Ki78,GS99,BP11}, the category $\RHB$ of oriented \dmnsnl{4} relative \hndlbds{2}  modulo \qvlnc{2}, is  identified with the category $\KTf$ of admissible Kirby tangles, representing the attaching maps of the 1- and 2-handles in $[0,1]^3$. Below we define the latter category and refer the reader to \cite{BP11, BBDP23} for all the details of the equivalence functor between $\KTf$ and $\RHB$. 

\begin{definition}
\label{kirby-admtangle/def}
For any integer $n \geqs 0$, let $E_n \subset [0,1]^2$  be a set of $2n$ points uniformly distributed along $\left] 0,1 \right[ \times \{ 1/2 \}$.  In particular, $E_0 = \emptyset$.
Then, given two integers $n,m \geqs 0$ an \textsl{admissible Kirby tangle} from $E_n$ to $E_m$ consists of the following data:
\begin{itemize}
\item[\(a)]
a collection of $k \geqs 0$ dotted unknots $U_1, U_2, \dots, U_k$, together with disjoint flat spanning disks $D_1, D_2, \dots, D_k$ embedded into $\left] 0,1 \right[^3$;
\item[\(b)]
an undotted tangle properly and smoothly embedded into $[0,1]^3$, which is transverse to the spanning disks, and whose components are endowed with the blackboard framing with respect to the projection $[0,1]^3 \to [0,1]^2$ that forgets the second coordinate. Moreover, the endpoints of each open component are two consecutive points in $(E_n \times \{0\}) $ or in $ (E_m \times \{1\})$.
\end{itemize}
Two admissible Kirby tangles are said to be {\sl 2-equivalent} if they are related by a 2-deformation, consisting of a finite sequence of ambient isotopies (including isotopies of the disjoint union of spanning disks of the dotted components) and the moves presented in Figure~\ref{Ktangles-hmoves/fig} which realise the corresponding handle moves. 

\begin{figure}[htb]
 \includegraphics{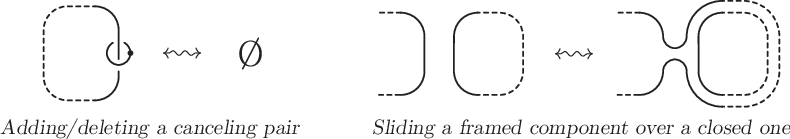}
 \caption{}
 \label{Ktangles-hmoves/fig}
\end{figure}
\end{definition}

\begin{figure}[htb]
 \centering
 \includegraphics{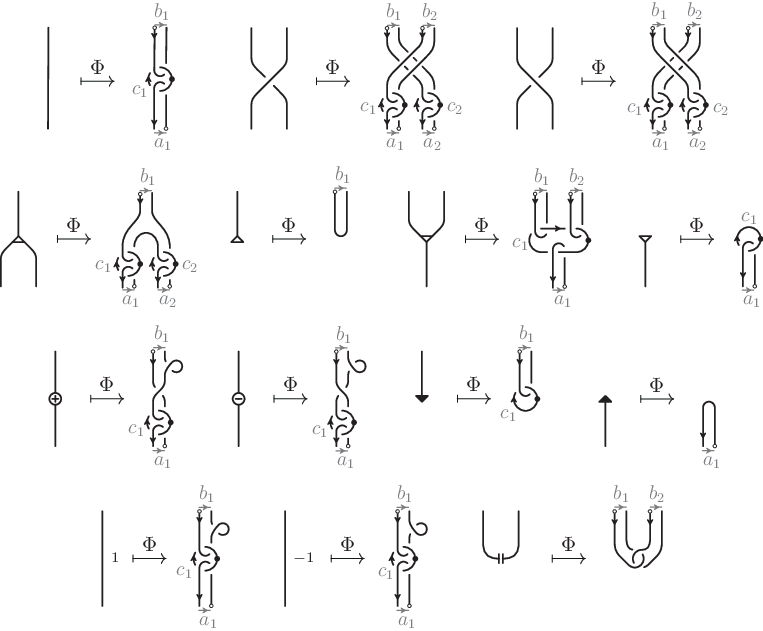}
 \caption{The functor $\Phi : \Algf \to \KTf$.}
 \label{phi01/fig}
\end{figure}

%\smallskip
2-equivalence classes of admissible Kirby tangles form a strict braided monoidal category $\KTf$ whose objects are the sets $E_{n}$ with $n \geqs 0$, and whose morphisms are \qvlnc{2} classes of admissible Kirby tangles. The composition  is given by vertical juxtaposition of the tangles, while the monoidal product is given by horizontal juxtaposition. The identity morphism on $E_1$ and the braiding morphisms $\gamma_{E_1,E_1}$, $\gamma^{-1}_{E_1,E_1}$ are the tangles in the first line in Figure~\ref{phi01/fig} (for now disregard the map $\Phi$, orientations and gray labels). Observe that we  present the Kirby tangles by plane projections which are injective on the spanning disks of the dotted components.

%\begin{figure}[htb]
% \includegraphics{Figures/KT-morph.eps}
% \caption{Identity and braiding morphisms in $\KTf$.}
% \label{KT-morph/fig}
%\end{figure}

Given a Kirby tangle, one can read a presentation of the fundamental group of the corresponding 2-handlebody by associating a generator to each oriented dotted component and a relator to each based oriented framed component \cite[Section 4.2]{GS99}. We will generalize this procedure to a functor from $\KTf$ to $\pres$.  

\begin{proposition}\label{ftop/thm}
There is a braided monoidal functor $\ftop: \KTf \to \pres$ which assigns to each morphism in $\KTf$ a relative presentation of the fundamental group of the corresponding 4-dimensional 2-handlebody.
\end{proposition}

\begin{proof}
Let $T:E_n\to E_m$ be a morphism in $\KTf$ with  $k \geqs 0$ dotted unknots $U_1, U_2, \dots, U_k$ spanning the disks $D_1, D_2, \dots, D_k$, and $s\geq 0$ framed components $L_1, L_2, \dots, L_s$. We choose an orientation of the spanning disks and a base point and orientation of the closed framed components. Each open framed component is oriented so that it runs from its left endpoint towards its right endpoint if such endpoints are in the target of $T$, and from its right endpoint towards its left endpoint if the endpoints are in the source. Moreover, the base point of such component is identified with the starting point  of the oriented arc. The above choices are illustrated in the tangle diagrams in Figure \ref{phi01/fig}. 

Let $\overline A=(a_1,\dots,a_n)$, $\overline B=(b_1,\dots,b_m)$, $C=\{c_1,\dots,c_k\}$. Observe that there is a bijective correspondence between $G=A\cup B\cup C$ and the set of generators of the fundamental group of the handlebody represented by $T$. Then to each framed component $L_i$, $1\leq i\leq s$ we associate a word $v_i$ in $C$ as follows: starting with $v_i=1$ and moving from the base point  in the positive direction, we add on the right the letter $c_j$ (resp. $c_j^{-1}$ ) each time $L_i$ passes through $D_j$ in the positive (resp. negative) direction. Let
$$
w_i=\begin{cases}
    v_i\;\text{ if $L_i$ is closed},\\
    a_k v_i \;\text{if $\partial L_i$  is the $k$th pair in the source},\\
    v_ib_k^{-1} \;\text{if $\partial L_i$  is the $k$th pair in the target.}
\end{cases} 
$$
Then we define
$$
\ftop(E_n)=n,\quad \ftop(T)=\langle\, (a_1,\dots,a_n),\,(b_1,\dots,b_m);\, \{c_1,\dots,c_k\}\, |\, w_1,\dots w_s\, \rangle.
$$
Different choices of orientations of the spanning disks and orientations and base points of the closed components, as well as changing $T$ by a 2-deformation (Figure \ref{Ktangles-hmoves/fig}) change $\ftop(T)$ by a sequence of AC-moves as follows:
\begin{itemize}
\item[-] reversing the orientation of $D_j$ applies \hrel{ac6};
\item[-] changing the base point of a closed component applies \hrel{ac4};
\item[-] reversing the orientation of a closed component applies \hrel{ac5};
\item[-] adding or deleting a canceling 1/2-pair applies \hrel{ac1};
\item[-]  2-handle slide replaces a relator by its product with a conjugate of another relator or its inverse, and is therefore implemented by a sequence of moves \hrel{ac2'}, \hrel{ac3}, \hrel{ac4} and \hrel{ac5};
\item[-] ambient isotopy changes the relators by a sequence of insertions/ cancellations of adjacent inverse letters and is implemented by \hrel{ac2'}.
\end{itemize}
To see that $\ftop$ preserves compositions, let $T:E_n\to E_m$ and $T':E_m\to E_{n'}$ and
\begin{gather*}
\ftop(T)=\langle\, (a_1,\dots,a_n),\,(h_1,\dots,h_m);\, c_1,\dots,c_k\, |\, w_1,\dots w_s\, \rangle\\
\ftop(T')=\langle\, (h_1,\dots,h_m),\,(b_1,\dots,b_{n'});\, c'_1,\dots,c'_{k'}\, |\, w'_1,\dots w'_{s'}\, \rangle.
\end{gather*} 
Then for each generator $h_i$, $1\leq i\leq m$, there exist $1\leq i_1\leq s$ and $1\leq i_2\leq s'$ such that in $\ftop(T')\circ\ftop(T)$, $h_i$ appears only in the words $w_{i_1}=v_{i_1}h_i^{-1}$ and $w_{i_2}'=h_iv_{i_2}'$, where  $v_{i_1}$ and $v_{i_2}'$ are words in the rest of the generators. Therefore, by applying \hrel{ac7} we eliminate $h_i$ and $w_{i_2}'$ and substitute $w_{i_1}$ with $v_{i_1}v_{i_2}'$. The resulting presentation is exactly $\ftop(T'\circ T)$. 

By definition $\ftop(T\otimes T')=\ftop(T)\otimes\ftop(T')$, hence it is left to verify that $\ftop$ preserves identities and braiding morphisms as follows:
\begin{gather*}
\ftop(\id_{E_1})=\langle\, (a_1),\,(b_1);\, c_1\, |\, a_1c_1^{-1},c_1b_1^{-1}\, \rangle=\langle\, (a_1),\,(b_1);\,\emptyset\, |\, a_1b_1^{-1}\, \rangle=\idp_1;\\
\ftop(\gamma_{E_1,E_1})=\langle\, (a_1,a_2),\,(b_1,b_2);\, c_1,c_2\, |\, a_1c_1^{-1},c_1b_2^{-1}, a_2c_2^{-1},c_2b_1^{-1}\, \rangle\\
=\langle\, (a_1,a_2),\,(b_1,b_2);\, \emptyset\, |\, a_1b_2^{-1}, a_2b_1^{-1}\, \rangle=\gammap_{1,1},
\end{gather*}
where once again we have used  \hrel{ac7} to eliminate the internal generators.
\end{proof}

%%%%%%%%%

The algebraic presentation of $\RHB$ is the category  $\Algf$, freely generated by a  BP Hopf algebra object. This algebraic structure was first defined and studied in \cite{BP11} in the general context of \textit{groupoid Hopf algebras}. The notion of a BP Hopf algebra was introduced in \cite{BD21, BBDP23} and corresponds to the special case of the trivial groupoid $\G = \{1\}$. 

\begin{table}[b]
  \includegraphics{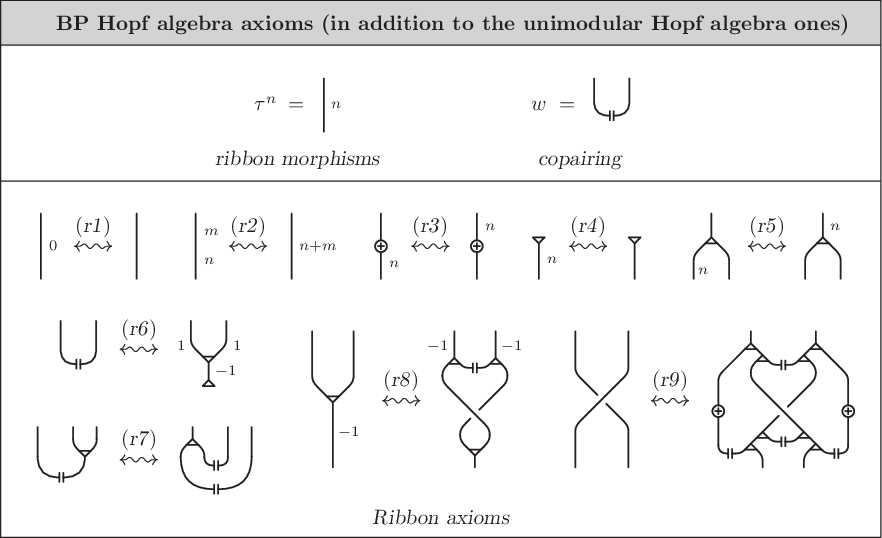}
\caption{}
\label{axiomsBPH/fig}\label{E:r1}\label{E:r2}\label{E:r3}\label{E:r4}\label{E:r5}\label{E:r6}\label{E:r7}\label{E:r8}\label{E:r9}\label{E:r1-2}
%\vskip-3pt
\end{table}

\begin{definition}\label{Hr/def}
Let $\calC$ be a braided monoidal category with tensor product $\otimes$, tensor unit $\one$, and braiding $\gamma$. Then a  \textit{BP Hopf algebra} in $\calC$ is a unimodular Hopf algebra $H$ equipped with the following structure morphisms (Table \ref{axiomsBPH/fig} where relations \hrel{r1-2} simply encode the notation for the powers of $\tau$):
\begin{itemize}
 \item a \textit{ribbon morphism} $\ribmorH : H \to H$ and its inverse $\ribmorH^{-1} : H \to H$;
 \item a \textit{Hopf copairing} $\copairH : \one \to H \otimes H $.
\end{itemize}
These structure morphisms are subject to the following axioms:
\begin{gather*}
 \antipH \circ \ribmorH = \ribmorH \circ \antipH, 
 \tag*{\hrel{r3}}
 \\
 \counH \circ \ribmorH = \counH, 
 \tag*{\hrel{r4}}
 \\
 \prodH \circ (\ribmorH \otimes \id) = \ribmorH \circ \prodH, 
 \tag*{\hrel{r5}} 
 \\
 \copairH = (\ribmorH \otimes \ribmorH) \circ \coprH \circ \ribmorH^{-1} \circ \unitH
 \tag*{\hrel{r6}}
 \\
 (\id \otimes \coprH) \circ \copairH = (\prodH \otimes \id_2) \circ (\id \otimes \copairH \otimes \id) \circ \copairH,
 \tag*{\hrel{r7}}
 \\
 \coprH \circ \ribmorH^{-1} =(\ribmorH^{-1} \otimes \ribmorH^{-1}) \circ \monH \circ \gamma^{-1} \circ \coprH,
 \tag*{\hrel{r8}}
 \\
 (\prodH \otimes \prodH) \circ (\antipH \otimes (\monH \circ \gamma^{-1} \circ \monH) \otimes \antipH) \circ (\rho_L \otimes \rho_R) = \gamma, 
 \tag*{\hrel{r9}}
\end{gather*}
where
\[
 \monH = (\prodH \otimes \prodH) \circ (\id \otimes \copairH \otimes \id) : H \otimes H \to H \otimes H
\]
is called the \textit{monodromy}\label{monodromy/def}, while the morphisms 
\[
 \rho_L = (\id \otimes \prodH) \circ (\copairH \otimes \id) : H \to H \otimes H 
 \quad \text{and} \quad
 \rho_R = (\prodH \otimes \id) \circ (\id \otimes \copairH) : H \to H \otimes H
 \]
define a left and a right $H$-comodule structure on $H$, respectively. 
\end{definition}

A graphical representation of the additional generators and defining axioms of a BP Hopf algebra can be found in Table~\ref{axiomsBPH/fig}, to be added to the list of generators and defining axioms of unimodular Hopf algebras given in Table~\ref{axiomsH/fig} (compare with \cite[Table~A1]{BBDP23}).

Let $\Algf$ be the braided monoidal category freely generated by a BP Hopf algebra $H$. 

\begin{remark}
    The axioms of unimodularity in Definition \ref{Hr/def} are slightly different, but equivalent to the ones in \cite[Definition 3.3.1]{BBDP23}. In particular, here we require  that the integral form and the integral element are two-sided (axioms \hrel{i1-1'} and \hrel{i2-2'}) which, together with the bialgebra axioms and \hrel{i3}, imply that they are $S$-invariant (the proofs of \hrel[i4]{i4-5} in Proposition \ref{antipode/thm} easily adapt to non-symmetric braiding). On the other hand, in \cite{BP11, BBDP23} it is assumed that the integrals are one-sided and $S$-invariant and it is proved that together with the bialgebra axiom and \hrel{i3} such relations imply that the integrals are two-sided. 
\end{remark}

We proceed with the proof of Theorem \ref{main2/thm}, which can be reformulated as follows:
\begin{theorem}\label{main21/thm}
There is a braided monoidal quotient functor $\falg:\Algf\to \calg$ such that the following diagram commutes:
\[
\begin{tikzcd}
\Algf \arrow{r}{\Phi} \arrow[swap]{d}{\falg} & \RHB \arrow{d}{\ftop} \\%
\calg \arrow{r}{\Psi}& \cwcob
\end{tikzcd}
\]
Moreover,  $\falg$ induces a braided monoidal category equivalence
$
\hat{\falg}:\Algf/\langle \tau=\id_H\rangle\to \calg,
$
where $\Algf/\langle \tau=\id_H\rangle$ is the quotient of $\Algf$ by the additional relation $\tau=\id_H$.
\end{theorem} 

\begin{proof}
It is proved in \cite[Theorem~4.7.5]{BP11}  that the  functor $\Phi: \Algf \to \KTf$ sending $H$ to $E_1$  and each structure morphism of $H$ to the corresponding Kirby tangle represented in Figure~\ref{phi01/fig} (disregarding tangle orientation and gray labels), is an equivalence of braided monoidal categories. A simpler  direct proof of this equivalence is given in \cite[Theorem~A]{BBDP23}.

We observe that any unimodular cocommutative  Hopf algebra in a symmetric category,  in particular $\calg$, is a BP Hopf algebra with trivial ribbon morphism and Hopf copairing ($\tau=\id_H$ and $\copairH=\eta\otimes\eta$). Indeed, in this case the additional ribbon axioms in Table \ref{axiomsBPH/fig} are automatically satisfied. Then, by the universal property of $\Algf$, there exists a functor $\falg:\Algf\to \calg$ sending the generating object of $\Algf$ to the generating object of $\calg$ and preserving structure morphisms.  

In order to show the commutativity of the diagram, it is enough to check that for each elementary morphism $F$ of $\Algf$ in Tables \ref{axiomsH/fig} and \ref{axiomsBPH/fig} we have $\ftop\circ \Phi(F)=\Psi\circ\falg(F)$. After recalling the definition of  $\ftop$ from the proof of Proposition \ref{ftop/thm} and the choices of orientation and base points presented in Figure \ref{phi01/fig}, this verification is straightforward and we only do it for $F=\mu,\Delta, \tau, \copairH$, leaving the rest to the reader:
   \begin{gather*}
       \ftop\circ\Phi (\mu)=\langle(a_1,a_2),(b_1);c_1,c_2\,|\, a_1c_1^{-1},\,a_2c_2^{-1},\,c_1c_2b_1^{-1}\rangle=
    \langle(a_1,a_2),(b_1);\emptyset\,|\, \,a_1a_2b_1^{-1}\rangle=\Psi\circ\falg(\mu); \\
    \ftop\circ\Phi (\Delta)=\langle(a_1),(b_1,b_2);c_1\,|\, a_1c_1^{-1},\,c_1b_1^{-1},\,c_1b_2^{-1}\rangle=
    \langle(a_1),(b_1,b_2);\emptyset\,|\, a_1b_1^{-1},\,a_1b_2^{-1}\rangle=\Psi\circ\falg(\Delta), \\
    \ftop\circ\Phi (\tau)=\langle(a_1),(b_1);c_1\,|\, a_1c_1^{-1},\,c_1b_1^{-1}\rangle=
    \langle(a_1),(b_1);\emptyset\,|\, a_1b_1^{-1}\rangle=\Psi(\id_1)=\Psi\circ\falg(\tau), \\
    \ftop\circ\Phi (\copairH)=\langle\emptyset,(b_1,b_2);\emptyset\,|\, b_1^{-1},\,b_2^{-1}\rangle=
    \langle\emptyset,(b_1,b_2);\emptyset\,|\, b_1,\,b_2\rangle=\Psi(\eta\otimes \eta)=\Psi\circ\falg(\copairH), 
   \end{gather*}
where in the first three lines line we use \hrel{ac7} to eliminate the internal generator(s), while in the last line we apply \hrel{ac5}.

Now we claim that  $\Algf/\langle \tau=\id_H\rangle$ is a symmetric category and $H$ is a unimodular cocommutative  Hopf algebra in it. Indeed, since $\tau=\id_H$, axioms \hrel[r3]{r3-5} in Table \ref{axiomsBPH/fig} become identities, while axioms \hrel{r6} in the same table and \hrel{a7} in Table \ref{axiomsH/fig}  imply that $\copairH =\eta\diam\eta$. Hence axiom \hrel{r7} becomes an identity as well, while axioms  \hrel{r8} and \hrel{r9} are reduced respectively to the cocommutative  axiom \hrel{a9} and to the symmetry axiom  in Table \ref{axioms2Alg/fig}. 

Hence  the universal property of $\calg$  gives a braided monoidal functor $\mathfrak{p}:\calg\to\Algf/\langle \tau=\id_H\rangle$ such that $\mathfrak{p}\circ\hat \falg$ and $\hat \falg\circ\mathfrak{p}$ fix the generating Hopf algebras and structure morphisms, and hence coincide with the corresponding identity functors.
 \end{proof}

%   \begin{figure}[htb]
%  \centering
%  \includegraphics{Figures/ftop1.eps}
%  \caption{Input data for $\ftop\circ\Phi (\mu)$ and $\ftop\circ\Phi (\Delta)$.}
%  \label{ftop/fig}
% \end{figure}

%%%%%%%%%%%%%%%%%%%%%%%%%%
\appendix

\section{The category $\calH$ freely generated by a cocommutative  Hopf algebra}\label{faith/sec}

Let $\calH$ be the symmetric monoidal category freely generated by a cocommutative  Hopf algebra $H$. In other words, $\calH$ has as set of objects $H^n$, $n\in \N$, and its morphisms are compositions of products of identities, $\Delta$, $\mu$, $\eta$, $\epsilon$ and $S$ modulo the bialgebra and antipode axioms in Table \ref{axiomsH/fig} (excluding the integral morphisms and axioms). Obviously, there is a symmetric monoidal  functor $\fh:\calH\to \calg$. 

 Pirashvili proved  that $\calH$ is equivalent to the opposite of the category of free groups $\Fopp$  \cite{Pi02}. A combinatorial proof of the `dual' statement, i.e. that the category of free groups $\Free$ is equivalent to the one freely generated by a commutative Hopf algebra, is given by Habiro in \cite{H16}.

We recall some facts about the category of free groups $\Free$ from \cite[Section 4]{H16}. $\Free$ is the full subcategory of the category of groups whose set of objects is $\{F_n\,|\, n\geq 0\}$. We will use the notation $F_n=\langle x_1\,\dots,x_n\rangle$ for the free group generated by the ordered set $(x_1\,\dots,x_n)$, which as an object of $\Free$,  will be identified with $n\in \N$. 
Then any morphism $\phi_w:\langle y_1\,\dots,y_m\rangle\to \langle x_1\,\dots,x_n\rangle$ in $\Free$ is uniquely specified by the sequence  of words $w=(w_1,\dots, w_m)$, $w_i=w_i(x_1,\dots, x_n)\in F_n$,  such that $\phi_w (y_i)=w_i$, $1\leq i\leq m$. Given another morphism $\phi_v:\langle z_1\,\dots,z_k\rangle\to \langle y_1\,\dots,y_m\rangle$, the composition $ \phi_w\circ\phi_v=\phi_{v(w)}$, where $(v(w))_j=v_j(w_1,\dots, w_m)$. 

 $\Free$ is a symmetric monoidal category whose monoidal structure is given by the free product of groups and the monoidal unit is $F_0$; in particular:
$$
n\otimes m=n+m,\quad\langle x_1\,\dots,x_n\rangle\otimes \langle y_1\,\dots,y_m\rangle=\langle x_1\,\dots,x_n,y_1\,\dots,y_m\rangle,
$$
and 
$$
\phi_{(w_1,\dots, w_n)}\otimes\phi_{(v_1,\dots, v_m)}=\phi_{(w_1,\dots, w_n, v_1,\dots, v_m)}.
$$
The symmetric structure is given by
$$
\gamma_{n,m}:\langle x_1\,\dots,x_n\rangle\otimes \langle y_1\,\dots,y_m\rangle \to \langle y_1'\,\dots,y_m'\rangle\otimes \langle x_1'\,\dots,x_n'\rangle=\phi_{(x_1'\,\dots,x_n', y_1'\,\dots,y_m')}.
$$

\begin{theorem}[Pirashvili]
There is a symmetric monoidal equivalence $ \Theta:\calH\to \Fopp$, which sends each elementary morphism $D:H^n\to H^m$ to $\phi_{w^D}:\langle y_1,...,y_m\rangle\to\langle x_1,...,x_n\rangle $, where
$$
w^\mu=(x_1x_2),\; w^\Delta=(x_1,x_1),\; w^\eta=1,\; w^\epsilon=\emptyset,\; w^S=(x_1^{-1}).
$$
\end{theorem}

The following proposition implies that the equivalence functor in Theorem \ref{main1/thm} can be interpreted as an extension of Pirashvili's functor to finite group presentations.

\begin{proposition}\label{free1/thm}
   There is a faithful functor $\ffree: \Fopp\to\pres$ such that the following diagram commutes. 
   \[
\begin{tikzcd}
\Fopp \arrow{r}{\ffree}  & \pres  \\%
\calH\arrow[swap]{u}{\Theta} \arrow{r}{\fh}& \calg \arrow{u}{\Psi}
\end{tikzcd}
\]
   In particular, the functor $\fh:\calH\to \calg$ is also faithful.
\end{proposition}

\begin{proof}
Define $\ffree$ to be the identity on  objects, while for any $\phi_w:\langle y_1\,\dots,y_m\rangle\to \langle x_1\,\dots,x_n\rangle$:
   $$
   \ffree(\phi_w)=\langle (x_1\,\dots,x_n), (y_1\,\dots,y_m);\emptyset\,|\, y_1^{-1}w_1,\dots y_m^{-1}w_m\rangle.$$
Observe that the map is well defined, since two words representing the same morphism in $\Free$ are related by insertions and cancellations of $gg^{-1}$ for some generator $g$, which induces \hrel{ac2'} on the image.  The fact that $\ffree$  respects the symmetric monoidal structure follows from its definition.  In order to see that it preserves compositions,  
let $\phi_w:\langle y_1\,\dots,y_m\rangle\to \langle x_1\,\dots,x_n\rangle$ and $\phi_v:\langle z_1\,\dots,z_k\rangle\to \langle y_1\,\dots,y_m\rangle$. Then
\begin{eqnarray*}
&\ffree(\phi_w\circ \phi_v)&=\langle (x_1\,\dots,x_n), (z_1\,\dots,z_k);y_1\,\dots,y_m\,|\, y_1^{-1}w_1,\dots y_m^{-1}w_m, z_1^{-1}v_1,\dots, z_k^{-1}v_k\rangle\\
&& =\langle (x_1\,\dots,x_n), (z_1\,\dots,z_k);\emptyset\,|\,  z_1^{-1}(v\circ w)_1,\dots, z_k^{-1}(v\circ w)_k\rangle,
\end{eqnarray*}
where we have used the first $m$ relators to eliminate the internal generators $y_i$, $1\leq i\leq m$. 

 Suppose now that $\phi_w, \phi_{w'}:\langle y_1\,\dots,y_m\rangle\to \langle x_1\,\dots,x_n\rangle$ and $\ffree(\phi_w)=\ffree(\phi_{w'}):n\to m$ in $\pres$.  
Let
$$
F_{n+m}=\langle x_1\,\dots,x_n,y_1\,\dots,y_m\rangle. 
$$
and $N_w$ and $N_{w'}$ be the normal subgroups of $F_{n+m}$ generated respectively by $y_1^{-1}w_1,\dots y_m^{-1}w_m$ and by $y_1^{-1}w'_1,\dots y_m^{-1}w'_m$. 
%%%%
 Since relative AC-moves induce isomorphisms of the presented groups preserving the images of all external generators and since both presentations have no internal generators, the AC-equivalence between $f(\phi_w)$ and $f(\phi_{w'})$ induces an isomorphism
$$ 
\alpha:F_{n+m}/N_w\to F_{n+m}/N_{w'}
$$
such that $\alpha\circ q_w=q_{w'}$, where $q_w,q_{w'}$ are the quotient maps. Therefore
$$ 
N_w=\ker q_w=\ker q_{w'}=N_{w'}. 
$$
On the other hand, since $f(\phi_w)$ and $f(\phi_{w'})$ are in the image of $\ffree$, both quotient groups are canonically identified with
$$ F_n=\langle x_1,\dots,x_n\rangle, 
$$
with the images of the $x_i$'s forming the standard free basis. Thus, from
$$ 
y_i^{-1}w_i=1=y_i^{-1}w_i' 
$$
in the common quotient, we obtain $w_i=w_i'$ in $F_n$. Hence $\phi_w=\phi_{w'}$.

We leave to the reader to check that $\Psi\circ\fh= \ffree\circ \Theta$, which implies the faithfulness of $\fh$.
\end{proof}

\end{document}